\documentclass{amsart}
\usepackage{amsmath,amsthm}
\usepackage{amsfonts,amssymb}

\newtheorem{thm}{Theorem}[section]

\newtheorem{lem}[thm]{Lemma}
\newtheorem{prop}[thm]{Proposition}

\theoremstyle{remark}
\newtheorem{rem}{Remark}[section]

\renewcommand{\theequation}{\thesection.\arabic{equation}}

\def\CB{{\mathcal B}}
\def\CF{{\mathcal F}}
\def\CJ{{\mathcal J}}
\def\CK{{\mathcal K}}
\def\CH{{\mathcal H}}
\def\CL{{\mathcal L}}
\def\CO{{\mathcal O}}
\def\CR{{\mathcal R}}
\def\CU{{\mathcal U}}
\def\CV{{\mathcal V}}
\def\CW{{\mathcal W}}

\def\a{{\mathbf a}}
\def\b{{\mathbf b}}
\def\m{{\mathbf m}}
\def\u{{\mathbf u}}
\def\v{{\mathbf v}}
\def\w{{\mathbf w}}
\def\x{{\mathbf x}}
\def\y{{\mathbf y}}
\def\z{{\mathbf z}}

\def\C{{\mathbb C}}
\def\H{{\mathbb H}}
\def\N{{\mathbb N}}
\def\R{{\mathbb R}}
\def\S{{\mathbb S}}
\def\W{{\mathbb W}}
\def\Z{{\mathbb Z}}

\def\eps{{\varepsilon}}

\def\1{\text{\bf {1}}}
\def\bs{\backslash}
\def\id{\mathop{\text{\rm{id}}}\nolimits}
\def\im{\mathop{\text{\rm{im}}}\nolimits}
\def\h{{\mathfrak h}}
\def\oline{\overline}
\def \la {\langle}
\def \ra {\rangle}

\newcommand{\wt}{\widetilde}
\newcommand{\wh}{\widehat}

\begin{document}

\title[Heat kernel transform]
{The heat kernel transform for the Heisenberg group}
\author{Bernhard  Kr\"otz, Sundaram Thangavelu and Yuan Xu}
\address{Department of Mathematics\\ University of Oregon\\
    Eugene, Oregon 97403-1222.}\email{kroetz@math.uoregon.edu}
\address{Stat-Math Division\\Indian Statistical Institute\\
8th Mile, Mysore Road\\Bangalore-560 059\\India.}
\email{veluma@isibang.ac.in}
\address{Department of Mathematics\\ University of Oregon\\
    Eugene, Oregon 97403-1222.}\email{yuan@math.uoregon.edu}

\date{\today}
\thanks{BK was supported in part by NSF grant DMS-0097314 and YX was
supported in part by NSF grant DMS-0201669. ST wishes to thank BK and
YX for the warm hospitality during his stay in Eugene.}

\begin{abstract}
The heat kernel transform $\CH_t$ is studied for the Heisenberg group in 
detail. The main result shows that the image of $\CH_t$ is a direct sum of
two weighted Bergman spaces, in contrast to the classical case of $\R^n$ 
and compact symmetric spaces, and the weight functions are found to be
(surprisingly) not non-negative. 
\end{abstract} 

\maketitle

\section{Introduction}
\setcounter{equation}{0}

Over the last decade one could observe interesting 
developments on the {\it heat kernel transform} 
for various 
types of homogeneous Riemannian manifolds $X$. 
Complete results have been obtained for 
compact Lie groups (cf.\ \cite{H, H2}) and, more
generally, for compact symmetric spaces
(cf.\ \cite{S}). 
For non-compact spaces $X$ the situation seems 
to be more complicated and little research 
has been undertaken in this direction: There 
is the well understood Euclidean case (e.g. $X=\R^n$, cf. \cite{B}) 
and some 
partial results have been obtained 
for non-compact Riemannian symmetric spaces (cf.\ \cite{KS}). 
The objective of this paper is to give a 
complete and self-contained discussion for the 
Heisenberg group.

\par Our concern is with the $(2n+1)$-dimensional
Heisenberg group $\H$ and its universal complexification 
$\H_\C$. For $t>0$ we write $k_t: \H \to\R^+$ 
for the heat kernel on $\H$. Contemplating on the spectral resolution of 
$k_t$, it is 
not hard to see that $k_t$ admits an analytic 
continuation  to a holomorphic function $k_t^\sim: 
\H_\C\to \C$. Consequently, for every 
$f\in L^2(\H)$ the convolution 
$f*k_t$ continues holomorphically to $H_\C$ and 
we obtain a map 
$$
\CH_t: L^2(\H)\to \CO(\H_\C), 
\ \ f\mapsto (f*k_t)^\sim\, .
$$
We refer to $\CH_t$ as the {\it heat 
kernel transform} on $\H$ with parameter $t>0$. 
The map $\CH_t$ is injective, left $\H$-equivariant
and becomes continuous if $\CO(\H_\C)$ is 
equipped with its natural Fr\'echet topology 
of compact convergence. It follows that 
$\im \CH_t$ is a reproducing kernel 
Hilbert space. Standard abstract arguments readily 
yield an  expression for the kernel function in terms of $k_{2t}^\sim$
(see (3.1.2) below). 
\par In all known cases (e.g. $X$ a 
compact symmetric space 
or $X=\R^n$) the image of the heat kernel transform has 
been a weighted Bergman space  $X_\C$ with regard 
to a postive weight function. 
It came 
to our surprise that the Heisenberg 
group deviates from this pattern. The main 
result of this paper asserts that 
\begin{equation}\label{1.1}
\im \CH_t =
\CB_t^+(\H_\C)\oplus \CB_t^-(\H_\C)
\end{equation}
is a direct sum of two weighted Bergman spaces on 
$H_\C$.  Most interestingly, 
the weight functions 
$W_t^{\pm}$ for $\CB_t^\pm (\H_\C)$ have 
an oscillatory nature and attain positive and 
{\it negative} values. This fact forces 
the use of a certain 
exhaustion  $\H_\C=\bigcup_{R>0} K_R$ to define 
the inner product 
a suitable dense subspace $\CV_t^\pm(\H_\C)$ of 
$\CB_t^\pm (\H_\C)$ by   
$$\la f, g\ra =\lim_{R\to \infty} \int_{K_R} f(z)
\oline{g(z)} W_t^\pm(z) \, dz\qquad 
(f,g\in 
\CV_t^\pm (\H_\C))\, ,$$
quite  reminiscent to the familiar notion of principal 
value. 

\medskip  Let us now describe the contents of this 
paper in more detail. In Section 2 we introduce 
our notation and recall some facts on the heat 
kernel $k_t$ on $\H$ and its analytic continuation 
to $\H_\C$. Subsequently in Section 3 we 
define the heat kernel transform and give a discussion 
of its general nature. 
\par For the remainder it is useful to identify 
$\H$ with $\R^{2n}\times \R$. In Section 4 we introduce 
for each spectral parameter $\lambda\in \R^\times$ a partial 
heat kernel transform 
$$H_t^\lambda : L^2(\R^{2n}) \to \CO(\C^{2n})$$
and show that $\im H_t^\lambda$ is a weighted 
Bergman space  $\CB_t^\lambda(\C^{2n})$ associated 
to an explicitly given positive weight function 
$W_t^\lambda\ : \C^{2n}\to \R^{2n}$. 
With these results we prove in Section 5 
that there is a natural left $\H$-equivariant equivalence 
\begin{equation}\label{1.2}
L^2(\H)\simeq \int_{\R^\times}^\oplus  
\CB_t^\lambda(\C^{2n})   
\, e^{2t\lambda^2} d\lambda\, .
\end{equation}
Moreover, within the identification (\ref{1.2}) the heat 
kernel transform $\CH_t$ becomes the diagonal 
operator $(H_t^\lambda)_\lambda$. 

\par In Section 6 we combine all previously obtained 
results to establish our  main result (\ref{1.1}). 
It turns out that the global weight functions 
$W_t^\pm$ admit an integral representation in terms 
of the partial weight functions $W_t^\lambda$. 
Finally, in the appendix we derive explicit expansions 
of  $W_t^\pm$ by Hermite polynomials and explain 
their oscillatory behavior. 
 
\medskip {\it Acknowledgement:} We would like to express our 
sincere gratitude to a referee who read the manuscript very 
carefully and pointed out several inaccuracies, gaps  and 
mistakes.

\renewcommand{\theequation}{\thesubsection.\arabic{equation}}

\section{The heat kernel on the Heisenberg group}
\setcounter{equation}{0}

\subsection{Notation}
\setcounter{equation}{0}

Let $\h$ denote the $(2n+1)$-dimensional
Heisenberg algebra with generators, say,
$$
X_1, \ldots, X_n, U_1, \ldots, U_n, Z
$$
and relations $[X_j,U_j]=Z$. 
In the sequel we will
often identify $\h$ with $\R^{2n+1}=\R^n\times \R^n
\times \R$. For that let $(\x,\u, \xi)$ with
$\x=(x_1, \ldots, x_n)$ and $\u=(u_1, \ldots, u_n)$
denote the canonical coordinates on $\R^{2n+1}$.
Then the map
$$
\R^{2n+1}\to \h, \ \ (\x,\u,\xi)\mapsto
\sum_{j=1}^n x_j X_j +\sum_{j=1}^n u_j U_j +\xi Z
$$
is a linear isomorphism providing us with
suitable coordinates for $\h$.

Let $\H$ denote a simply connected Lie group
with Lie algebra $\h$, the {\it Heisenberg group}.
We will identify $\H$ with $\h$ through
the exponential function $\exp=\id: \h\to \H$.
As $\H$ is two step,
the Baker-Campbell-Hausdorff formula provides
the group law
$$
(\x,\u,\xi)(\x',\u',\xi')= (\x+\x', \u+\u', {\frac12}(\x\cdot \u'-
\u\cdot \x') + \xi+\xi')\, .
$$
Here $\x\cdot \u=\sum_{j=1}^n x_j u_j$, as usual,
denotes the standard pairing on $\R^n$.
We notice in particular that
$$
(\x,\u,\xi)^{-1}=(-\x,-\u,-\xi)\, .\leqno(2.1.1)
$$

Write $dh$ for a Haar measure on $\H$. We can and
will normalize $dh$ in such a way that it
coincides with the product of Lebesgue
measures, i.e.
$$
\int_\H f(h) \,dh =\int_{\R^{2n+1}} f(\x,\u,\xi)
\,d\x \,d\u \,d\xi
$$
for all $f\in C_c(\H)$.

Write $\H_\C$ for the universal
complexification of $\H$. Of course we can
identify $\H_\C$ with $\C^{2n+1}$ and we will
often do so.  We will write $(\z,\w,\zeta)$
for the coordinates on $\C^{2n+1}$ where
$\z=\x+i\y$, $\w=\u+i\v$ and $\zeta=\xi+i\eta$.

For any simply
connected nilpotent Lie group $H$ the
polar mapping
$$
H\times \h\to H_\C, \ \ (h,X)\mapsto h\exp(iX)
$$
is a homeomorphism. Furthermore the Haar measure on
$H_\C$ decomposes as
$$
\int_{H_\C} f(g)\,dg =\int_H\int_\h
f(h\exp(iX))\,dX \,dh \leqno(2.1.2)
$$
for all $f\in C_c(H_\C)$.

For the Heisenberg group $\H$, the polar mapping
is explicitly given by
$$
\left( (\x, \u, \xi),(\x',\u',\xi')\right)
\mapsto (\x+i\x', \u+i\u', \frac{i}{2} (\x\cdot \u' -
\u\cdot \x') +\xi +i\xi')
$$
where $h = (\x,\u,\xi)$ and $X = (\x',\u',\xi')$.
In particular the Haar measure on $\H_\C$
can be chosen as the product of Lebesgue measures
$d\x\,d\y\,d\u\,d\v\,d\xi\,d\eta$.

For integrable functions $f, g$ on $\H$ we
define their convolution by
$$
(f*g)(x)=\int_\H f(h)g(h^{-1}x) \,dh \qquad (x\in \H)\, .
$$
In coordinates this is explicitly given by
$$
(f*g)(\x,\u,\xi)=\int_{\R^{2n+1}}
f(\x',\u',\xi') g\left((-\x',-\u',-\xi')(\x,\u,\xi)\right)
\,d\x'\,d\u'\,d\xi'\, .
$$

\subsection{The heat kernel}
\setcounter{equation}{0}

Write ${\CU}(\h)$ for the universal
enveloping algebra of $\h$ and define the
Laplace element in ${\CU}(\h)$ by
$$
{\CL}=\sum_{j=1}^n X_j^2 +\sum_{j=1}^n U_j^2 + Z^2\, .
$$

For $X\in \h$ we write $\tilde X$ for the
left invariant vector field on $\H$, i.e.,
$$
(\tilde X f)(h)=\frac{d}{dt}\Big|_{t=0} f(h\exp(tX))
$$
for $f$ a function on $\H$ which is differentiable
at $h\in \H$. Write $\rho$ for the right regular
representation of $\H$ on $L^2(\H)$, i.e.
$$
(\rho(h)f)(x)=f(xh)
$$
for $h,x\in \H$ and $f\in L^2(\H)$. With
$d\rho$ the derived representation we then have
$d\rho(X)=\tilde X$ for all $X\in \h$. In particular
if
$$
\Delta=\sum_{j=1}^n \tilde X_j^2 +\sum_{j=1}^n \tilde
U_j^2 +\tilde Z^2
$$
denotes the Laplace operator on $\H$, then
$d\rho({\CL})=\Delta$.  

Set $\R^+=(0,\infty)$. Our concern will be with the heat equation
on $\H\times \R^+$
$$
\partial_t u(h,t)=\Delta u(h,t)
$$
for appropriate functions $u(h,t)$ on $\H\times \R^+$.
The fundamental solution is
given by the heat kernel $k_t(h)$ which can
be computed as follows:
$$
k_t(\x,\u,\xi)=c_n  
\int_\R e^{-i\lambda\xi} e^{-t\lambda^2}
\left( \frac{\lambda}{\sinh \lambda t}\right)^n
e^{-\frac{\1}{4} \lambda (\coth t\lambda)
(\x\cdot \x+ \u\cdot \u)} \,d\lambda \leqno(2.2.1)
$$
with $c_n=(4\pi)^{-n}$ (this follows from
a slight modification of \cite[Theorem 2.8.1]{T}.)
It satisfies the usual property of
$k_t * k_t = k_{2t}$ (see, for example, \cite[(2.87) and Corollary 2.3.4]{T}).

If $f$ is an analytic function on $\H$
which holomorphically extends to $\H_\C$, then we
write $f^\sim$ for this holomorphic extension.
The explicit formula (2.2.1) now implies
that $k_t$ has a holomorphic continuation to
$\H_\C$ which is given by
$$
k_t^\sim(\z,\w,\zeta)=c_n  
\int_\R e^{-i\lambda\zeta} e^{-t\lambda^2}
\left( \frac{\lambda}{\sinh \lambda t}\right)^n
e^{-\frac{\1}{4} \lambda (\coth t\lambda)
(\z\cdot \z+ \w\cdot \w)} \,d\lambda \leqno(2.2.2)
$$
for $(\z,\w,\zeta)\in \C^{2n+1}=\H_\C$.
It follows from (2.1.1) and (2.2.2) that
$$
k_t^\sim (z)=k_t^\sim (z^{-1}) \qquad (z\in \H_\C)\, .
\leqno(2.2.3)
$$
Furthermore, as $k_t\geq 0$ is real,
we record
$$
\oline{k_t^\sim (z)}=k_t^\sim(\oline z)
\qquad (z\in \H_\C)\, \leqno(2.2.4)
$$
Here, as usual,  $z\mapsto \oline z$ denotes the complex conjugation 
of $\H_\C$ with respect to the real form $\H$.

\section{The heat kernel transform}
\setcounter{equation}{0}

\subsection{Definition and basic properties}
\setcounter{equation}{0}

Let $C\subseteq \H_\C$ be a compact subset.
Then it follows from (2.2.2) that
$$
\sup_{z\in C} \int_{\H} |k_t^\sim(h^{-1}z)|^2
\,dh <\infty\, .\leqno(3.1.1)
$$
Fix $t>0$. Then (3.1.1) implies  
that $f*k_t$ has an analytic continuation to
$\H_\C$ for all $f\in L^2(\H)$.
In particular  we obtain a
linear map
$$
\CH_t: L^2(\H)\to {\CO}(\H_\C), \ \
f\mapsto (f*k_t)^\sim; \ \CH_t(f)(z)=\int_{\H_\C}
f(h)k_t^\sim (h^{-1} z) \,dh \  .
$$
We will call $\CH_t$ the {\it heat kernel transform}.

In the sequel we wish to consider ${\CO}(\H_\C)$ as a
Fr\'echet space -- the topology being the one of compact
convergence. If $h\in H$ and $f$ is a function on $\H$ or $\H_\C$,
then we write $\tau(h)f=f(h^{-1}\cdot)$.
The following properties of $\CH_t$  are  immediate:

\begin{itemize}
\item $\CH_t$ is continuous (because of (3.1.1))
\item  $\CH_t$ is injective (note that
$\CH_t(f)=e^{t\Delta}f$ and $\Delta$ is a negative
definite operator).
\item $\CH_t$ is $\H$-equivariant, i.e.
$\CH_t\circ \tau(h)= \tau(h) \circ \CH_t$ for all
$h\in \H$ (this is a general fact for the convolution
on a locally compact group).
\end{itemize}

We will endow $\im \CH_t$ 
with the Hilbert topology induced from $L^2(\H)$.
As $\CH_t$ is continuous we see that $\im \CH_t$
is an $\H$-invariant Hilbert space of holomorphic
functions on $\H_\C$. As such $\im \CH_t$ has
continuous point evaluations, i.e. for all
$z\in \H_\C$ the map
$$
\im \CH_t\to \C, \ \ f\mapsto f(z)
$$
is continuous. Hence $f(z)=\la f, {\CK}_z^t\ra$
for a unique element  ${\CK}_z^t\in \im \CH_t$.
We then obtain a positive definite kernel function
$$
{\CK}^t: \H_\C\times \H_\C\to \C;\  
{\CK}^t(z,w)=\la {\CK}_w^t, {\CK}_z^t\ra={\CK}_w^t(z)
$$
which is holomorphic in the first and anti-holomorphic
in the second variable. Moreover, the $\H$-invariance
of $\im \CH_t$ translates into  
${\CK}^t(hz,hw)={\CK}^t(z,w)$ for all $h\in \H$
and $z,w\in \H_\C$.

Let us compute ${\CK}^t$. Fix $w\in \H_\C$.
Let $g\in \im \CH_t$. Then $g=\CH_t(f)$
for some $f\in L^2(\H)$ and
$$
\la g, {\CK}_w^t\ra =g(w)=\CH_t(f)(w)
=(f*k_t)^\sim(w)=\int_\H f(h)k_t^\sim (h^{-1}w)\,dh\, .
$$
As this holds for all $g\in \im \CH_t$, we thus
conclude that
$$
\CH_t^{-1}({\CK}_w^t)(h)=\oline {k_t^\sim(h^{-1}w)}
=k_t^\sim(\oline w^{-1} h) \qquad (h\in \H)
$$
where for the last equality we used the facts (2.2.3-4).
From this we now get for all $w,z\in \H_\C$ that
\begin{align*}
\CK_w^t(z) & = \CH_t(k_t^\sim(\oline w^{-1}\cdot) )(z)=
\int_\H k_t^\sim (\oline w^{-1}h) k_t^\sim(h^{-1}z)\, dh \\
& =\int_\H k_t(h) k_t^\sim(h^{-1}\oline w^{-1} z)\, dh \\
&=(k_t*k_t)^\sim (\oline w^{-1} z) \\
&=k_{2t}^\sim(\oline w^{-1}z).
\end{align*}
We have thus shown that the kernel function is
given by
$$
{\CK}^t(z,w)=k_{2t}^\sim (\oline w^{-1}z)
\qquad (z,w\in \H_\C)\, .
\leqno (3.1.2)
$$

\subsection{General remarks on integral transforms and Bergman spaces}
\setcounter{equation}{0}

The setup for this Section is as follows: We let 
$N$ be a positive integer and $G$ be a Lie group 
which acts on $\R^N$ in a measure preserving manner. 
We assume that the action of $G$ extends to an action 
on $\C^N$ by measure preserving biholomorphisms. 
Our next data is 
a continuous (integral) transform 
$$\Phi: L^2(\R^N)\hookrightarrow \CO(\C^N)$$
which we assume to be $G$-equivariant. In this 
way $\im \Phi$ becomes a $G$-invariant Hilbert space 
of holomorphic functions on $\C^N$. We write 
$\CK: \C^N\times \C^N\to \C$ for the 
corresponding kernel function. 

\medskip\noindent
\textbf{Example 3.1.}
(a)  The heat kernel transform 
$\CH_t: L^2(\H)\to \CO(\H_\C)$ meets the general assumptions
from above. In fact, for  
$N=2n+1$ we may identify $\H$ with $\R^N$ and 
$\H_\C$ with $\C^N$. Furthermore the group $G=\H$
acts from the left on $\H=\R^N$ and $\H_\C=\C^N$ in a 
measure preserving manner. 
\par\noindent  (b) The partial heat 
kernel transforms $H_t^\lambda: L^2(\R^{2n})\hookrightarrow
\CO(\C^{2n})$ introduced in Section 4 below satisfy 
the general assumptions made above. 
\qed

\medskip

For the remainder we will assume that
$\im \Phi={\CB}(\C^N, W)$ is
a weighted Bergman space
for some measurable weight function
$W: \C^N \to \R$, i.e.,
$$\CB(\C^N, W)=\{ f\in {\CO}(\C^N): 
\int_{\C^N} |f(z)|^2 \ |W(z)|\, dz <\infty\}
$$
Hilbert structure given by  
$$\la f, g\ra =\int_{\C^n} f(z) \oline{g(z)} W(z) dz\leqno(3.2.1)$$

As the action of $G$ on
$\CB(\C^N,W) $ is unitary, the weight function
$W$ should be left $G$-invariant, i.e.
$$
W(g.z)=W(z) \qquad (g\in G, z\in \C^N)\, .
\leqno(3.2.2)
$$

\par What we cannot expect however is that $W$ is non-negative. 
It might then be a surprise
that (3.2.1) still defines a Hilbert structure. 
As the following  example shows,
this is a phenomenon which already appears in one
variable.

\medskip\noindent
\textbf{Example 3.2.}
We consider the unit disk
$D=\{ z\in \C : |z|<1\}$. For a measurable subset $A\subseteq D$ 
write $\1_A$ for its characteristic function. Define a weight function $W$
on $D$ by 
$$ 
W=\1_{\{\frac12 \leq |z|<1\}}- \1_{\{|z|<\frac12\}}\, . 
$$
With $W$ we form the weighted Bergman space
$$
\CB^2(D,W):=\{ f\in \CO(D): \int_D |f(z)|^2 
      \ |W(z)| \,dx\,dy <\infty\}
$$
and endow it  with the sesquilinear bracket 
$$\la f, g\ra =\int_D f(z) \oline {g(z)}\  W(z)\,  dx\, dy\, . $$

We will show that  $(\CB^2(D,W), \la\cdot, \cdot\ra)$  is a Hilbert space. For that
we first observe that $\{z^n\}_{n\in \N_0}$ is an orthogonal
system in $\CB^2(D,W)$. This is because $W$ is rotationally
invariant. Next we compute 
\begin{align*}
\la z^n, z^n \ra &=2\pi \int_{\frac12}^ 1 r^{2n+1} \,dr - 2\pi \int_0^{\frac12}
r^{2n+1} \,dr\\
&=\frac{\pi}{n+1} \left[1- \left(\frac12\right)^{2n+1}\right]>0
\end{align*}  
for any $n\in \N_0$. Thus if $f=\sum_n a_n z^n\in \CB^2(D,W)$ is an
arbitrary element, then
\begin{equation}\label{zix}
\la f, f\ra=\sum_n |a_n|^2 \frac{\pi}{n+1} \left[1- \left(\frac12\right)^{2n+1}
\right] \geq 0\end{equation}
and $\la f, f\ra =0$ if and only if $f=0$. 
This shows that $\la\cdot, \cdot\ra$ defines a pre Hilbert structure 
on $\CB^2(D,W)$. 
Next notice that 
\begin{equation}\label{zax}
\int_D |f(z)|^2\ |W(z)| \,dx\,dy =\sum_n |a_n|^2 \frac{\pi}{n+1}\, .\end{equation}
It follows from  identiies  (\ref{zix}) and (\ref{zax}) that 
$\la \cdot, \cdot\ra$ and 
the Hilbert bracket 
$(f|g)=\int_D f(z) \oline {g(z)}\  |W(z)|\,  dx\, dy$ induce equivalent 
norms. Hence $(\CB^2(D,W), \la\cdot, \cdot\ra)$ is 
a Hilbert space. 
\par  Finally we note that $W$ is uniquely characterized by the Hilbert norm
on $\CB^2(D,W)$, i.e. $\CB^2(D,W)=\CB^2(D,W')$ if and only if $W=W'$ almost 
everywhere (use Stone-Weierstra\ss).
\qed

\medskip

We conclude this section with some general remarks on how to obtain 
the weight function $W$. Define a subspace of $\im \Phi$ by 
$$
(\im \Phi)_0=\text{\rm span}\{ \CK_x: x\in \R^N\}\, .
$$
Since a holomorphic function on $\C^N$ which vanishes 
on $\R^N$ is identically zero, we conclude that 
$(\im \Phi)_0$ is dense in $\im \Phi$. Hence 
$\im \Phi =\CB(\C^N, W)$ will hold precisely if 
$$
\CK(x, x')=\la \CK_{x'}, \CK_x\ra =\int_{\C^N}
\CK_{x'}(z) \oline{\CK_x(z)} W(z)\,dz \leqno(3.2.3)
$$  
for all $x,x'\in \R^N$.  The formula (3.2.3) is 
actually quite helpful and will be applied  in Section 4 below.

\section{The $\lambda$-twisted heat-kernel transform}
\setcounter{equation}{0}

For $\lambda\in \R$, $\lambda\neq 0$,  we will  introduce 
a $\lambda$-twisted heat kernel transform 
$H_t^\lambda: L^2(\R^{2n})\to \CO(\C^{2n})$.
We will show that the image of $H_t^\lambda$ is 
a weighted Bergman space $\CB_t^\lambda (\C^{2n})$ on $\C^{2n}$. 
Further we provide an inversion formula for $H_t^\lambda$. 
\par The results of this section are the building blocks 
for our general discussion of the heat kernel transform 
$\CH_t: L^2(\H)\to \CO(\H_\C)$ in the following
sections.

\subsection{Notation}
\setcounter{equation}{0}

Let $\lambda\in \R$, $\lambda\neq 0$. For suitable functions $F$ on $\H$ we 
define a function $F^\lambda$ on $\R^{2n}$ by  
$$
F^\lambda(\x,\u)=\int_\R e^{i\lambda\xi}  F(\x,\u,\xi)\, d\xi\, .
$$

\par For $f, g\in L^1(\R^{2n})$ the {\it $\lambda$-twisted convolution} is 
defined by 
$$
 (f *_\lambda g) (\x,\u) = \int_{\R^{2n}} f(\x',\u') g(\x-\x',\u-\u')
        e^{-i \frac{\lambda}{2} (\x' \cdot \u-\x \cdot \u')}\,d\x'\,d \u'. 
$$
Notice that we have for Schwartz functions 
$F,G\in S(\H)=S(\R^{2n+1})$ that 
\begin{equation}\label{eins}  
(F*G)^\lambda =F^\lambda *_\lambda G^\lambda\, .
\end{equation} 

\par Let $\Delta_{\text {\rm sub}}=d\rho\left( \sum_{j=1}^n 
(\tilde X_j^2+\tilde Y_j^2)\right)$
denote the sublaplacian on $\H$. The heat kernel for
$\Delta_{\text {\rm sub}}$ 
is denoted by $p_t$ and its inverse Fourier 
transform in the central variable is explicitly given by 
\begin{equation}\label{eq:4.1}
p_t^\lambda(\x,\u) = c_n \left(\frac{\lambda}{\sinh t\lambda}\right)^n
e^{-\frac{\lambda}{4}\coth(\lambda t)(|\x|^2+|\u|^2)}
\end{equation}
with $c_n = (4\pi)^{-n}$. 

\par For all $f\in L^2(\R^{2n})$ the twisted convolution 
$f *_\lambda p_t^\lambda$ 
has an analytic continuation to $\C^{2n}$. 
In particular, there is a {\it $\lambda$-twisted heat kernel transform}
$$ 
H_t^\lambda: L^2(\R^{2n}) \to \CO (\C^{2n}), \ \ f \mapsto 
    (f *_\lambda p_t^\lambda)^\sim\, .    
$$
In coordinates we have 
$$
H_t^\lambda (f)(\z,\w)= \int_{\R^{2n}} f(\x',\u')p_t^\lambda(\z-\x',\w-\u') 
     e^{-\frac{i}{2} \lambda(\x'\cdot \w-\u'\cdot \z )}\,d\x'\,d\u'\, . 
$$
We define a unitary representation $\tau^\lambda$ of $\R^{2n}$ on 
$L^2(\R^{2n})$ by 
$$
(\tau^\lambda(\a,\b) f)(\x,\u)=
e^{-\frac{i\lambda}{2} (\a\cdot \u -\b\cdot \x)}f(\x-\a,\u-\b)
$$
for $(\a,\b)\in \R^{2n}$, $f\in L^2(\R^{2n})$ and $(\x,\u)\in \R^{2n}$. 
Likewise $\tau^\lambda$ defines an action of $\R^{2n}$ on $\CO(\C^{2n})$ 
via 
$$
(\tau^\lambda(\a,\b) f)(\z,\w)=e^{-\frac{i\lambda}{2}(\a\cdot \w -\b\cdot \z)}
 f(\z-\a,\w-\b)
$$
where  $(\a,\b)\in \R^{2n}$, $f\in \CO(\C^{2n})$ and $(\z,\w)\in \C^{2n}$. 

\par As for functions $F, G\in L^1(\H)$ we have $\tau(h)F* G=\tau(h)(F*G)$ 
for all $h\in \H$, it is immediate from (\ref{eins}) that  $H_t^\lambda$ 
becomes $\R^{2n}$-equivariant, i.e 
\begin{equation} \label{eq:a}
H_t^\lambda (\tau^\lambda(\a, \b) f)= \tau^\lambda(\a, \b) (H_t^\lambda (f))
\end{equation}
for all $(\a, \b)\in \R^{2n}$ and $f\in L^2(\R^{2n})$. 


\begin{rem} For the proofs in the sequel it is notationally 
convenient to prove the assertions for the ``essential case'' $\lambda=1$ 
only. Whenever  
we do so we will use a simplified notation: we write 
$f\times g$ instead of $f*_1 g$ for the $1$-twisted 
convolution; further we will drop all sub- and superscripts 
involving $\lambda=1$, i.e. $p_t^1$ becomes $p_t$, $H_t^1$
becomes $H_t$ etc. 
\end{rem}

\subsection{Determination of the weight function}
\setcounter{equation}{0}

Our objective is to find a non-negative  weight 
function $W_t^\lambda$ on $\C^{2n}$ such that 
\begin{align} \label{eq:4.2}
 \int_{\C^{2n}} |H_t^\lambda (f)(\z,\w) |^2  W_t^\lambda(\z,\w) \,d\z \,d\w
=\int_{\R^{2n}} |f(\x,\u)|^2 \,d\x\,d\u \end{align}
for all $ f\in L^2(\R^{2n}).$

\begin{prop} \label{prop:4.1}
A weight function $W_t^\lambda $ which satisfies  \eqref{eq:4.2} is
given by
\begin{equation} \label{eq:4.3}
W_t^\lambda(\x+i\y,\u+i\v) =4^n e^{\lambda(\u \cdot \y-\v \cdot \x)} 
p_{2t}^\lambda(2\y,2\v).
\end{equation} 
\end{prop}

\begin{rem} The weight function  $W_t^\lambda$ is unique in the sense that 
is the unique measurable function $W_t^\lambda: \C^{2n}\to \R_{\geq 0}$ which satisfies
(\ref{eq:4.2}). This will be shown in Lemma \ref{lem=unique} below. 
\end{rem}

\begin{proof} 
We restrict our attention to the case $\lambda=1$. As mentioned earlier 
we will write now $p_t$ and $W_t$ in place of $p_t^\lambda$ and $W_t^\lambda$, 
respectively, 
and write $ f\times g $ for the $1-$twisted convolution of $f$ and $g$. 
Via $H_t$ we can transfer the Hilbert space structure of 
$ L^2( \R^{2n} )$ to $\im H_t$ and make it into  
Hilbert space of holomorphic functions. 
Write $ K^t(\z,\w;\z',\w') $ for the corresponding reproducing kernel. 
Arguing as in Subsection 3.2, the inner product $\la\cdot,\cdot\ra_t $ on the 
image is uniquely determined by the equality
\begin{equation}\label{zwei}
K^t(\a,\b;\a',\b') = \la K^t_{(\a,\b)},K^t_{(\a',\b')}\ra_t
\end{equation}
for all real pairs $ (\a,\b), (\a',\b')\in \R^n\times \R^n$. 

As the heat kernel transform $ f \rightarrow H_t(f)= (f\times p_t)^\sim 
$ commutes with
the twisted translation (see equation (\ref{eq:a})),  
we may assume $ (\a',\b') = 0 $ in (\ref{zwei}). 
As $p_t\times p_t=p_{2t}$, arguing as in Subsection 3.1 
readily yields 
$$
K^t_{(\a,\b)}(\z,\w) = p_{2t}(\z-\a,\w-\b) 
      e^{-\frac{i}{2}(\a\cdot\w-\b\cdot\z)}\, .
$$
In particular, $K^t_{(0,0)}=p_{2t}$ and $K^t(\a,\b, 0, 0)=p_{2t}(\a,\b)$. 
Thus (\ref{zwei}) translates into 
$$
p_{2t}(\a,\b) = \int_{\C^n} \int_{\C^n}
p_{2t}(\z-\a,\w-\b) e^{-\frac{i}{2}(\a\cdot\w-\b\cdot\z)}
\overline{p_{2t}(\z,\w)} W_t(\z,\w) \,d\z \,d\w\, .
$$
This is established in Lemma \ref{lem:4.2} below.
\end{proof}

\begin{lem} \label{lem:4.2}
For $\a, \b \in \R^n $ we have
$$ \int_{\C^n}\int_{\C^n} p_{2t}^\lambda(\z+\a,\w+\b) 
e^{\frac{i\lambda}{2}(\a\cdot\w-\b\cdot\z)}
\overline{p_{2t}^\lambda(\z,\w)} W_t^\lambda(\z,\w) \,d\z\,d\w
 = p_{2t}^\lambda(\a,\b).
$$
\end{lem}

\begin{proof}  We will prove the assertion for $\lambda=1$. 
Further, by  the product nature of the functions involved, we may assume in 
addition that $n = 1$. 
\par Expanding out and simplifying we have
\begin{align*}
& p_{2t}(x+a+iy,u+b+iv) \overline{p_{2t}(x+iy,u+iv)}
= (4\pi)^{-2} (\sinh 2t)^{-2} e^{-\frac{1}{2}(\coth 2t) (x^2+u^2)} \\
&\qquad \qquad\qquad
\cdot  e^{-\frac{1}{4}(\coth 2t) (a^2+b^2)}e^{\frac{1}{2}(\coth 2t) (y^2+v^2)}
e^{-\frac{1}{2} (\coth 2t) \left(
a(x+iy)+b(u+iv) \right)}.
\end{align*}
We can combine the terms $e^{-\frac{1}{2}(\coth 2t) (x^2+u^2)} $ and
$$
e^{(uy-vx)} = e^{(\coth 2t) (uy \tanh (2t)  -xv \tanh (2t))}
$$
to get
\begin{align*}
& p_{2t}(x+a+iy,u+b+iv) \overline{p_{2t}(x+iy,u+iv)}
e^{(uy-vx)}\\
& \hspace{.5in} = (4\pi)^{-2} 
(\sinh 2t)^{-2} e^{-\frac{1}{4}(\coth 2t) (a^2+b^2)}
e^{\frac{1}{2}(\coth 2t +\tanh 2t)(y^2+v^2)} \\
& \hspace{.7in} \cdot e^{-\frac{1}{2}(\coth 2t) \left((x+v\tanh (2t) )^2+
(u-y\tanh (2t)  )^2\right)} e^{-\frac{1}{2} (\coth 2t) \left(
a(x+iy)+b(u+iv)\right)}.
\end{align*}
Using the identity $ \tanh 2t + \coth 2t = 2 \coth 4t $ and simplifying further
we get
\begin{align*}
& p_{2t}(z+a,w+b) e^{\frac{i}{2}(aw-bz)}\overline
     {p_{2t}(z,w)} W_t(z,w)\\
&\hspace{.3in} =  4^{-2} \pi^{-3}
 (\sinh 2t)^{-3} e^{-\frac{1}{8}\coth 2t (a^2+b^2)}
e^{(\coth 4t -\coth 2t)(y^2+v^2)} \\
& \hspace{.5in}\cdot
e^{-\frac{1}{2}(\coth 2t) \left((x+ \frac{a}{2}+\tanh (2t) v)^2+
(u+\frac{b}{2}-y\tanh (2t) )^2\right)} e^{-\frac{i}{2} (\coth 2t) (ay+bv)}
e^{\frac{i}{2}(au-bx)},
\end{align*}
where $z = x+iy$ and $w=u+iv$. 

First consider the integral
\begin{align*}
&\int_{\R^2} e^{\frac{i}{2}(au-bx)}e^{-\frac{1}{2}(\coth 2t)
\left((x+ \frac{a}{2}+v \tanh (2t) )^2+ (u+\frac{b}{2}- y \tanh (2t))^2\right)}
\,dx\,du\\
&\qquad
= e^{\frac{i}{2} (\tanh 2t) (ay+bv)} \int_{\R^2} e^{\frac{i}{2}(au-bx)}
e^{-\frac{1}{2}(\coth 2t) (x^2+u^2)}\,dx\,du\\
&\qquad = 2 \pi (\tanh 2t) e^{\frac{i}{2} (\tanh 2t) (ay+bv)}
  e^{-\frac{1}{8}(\tanh 2t) (a^2+b^2)}.
\end{align*}
Up to an explicit factor the remaining integral is
$$
\int_{\R^2} e^{-\frac{i}{2}(\coth 2t-\tanh 2t)(ay+bv)}
e^{-(\coth 2t -\coth 4t)(y^2+v^2)} \,dy\,dv.
$$
As $\coth 2t -\tanh 2t = 2(\sinh 4t)^{-1} $ and $ \coth 2t -\coth 4t =
(\sinh 4t)^{-1} $ the above integral reduces to
$$
\int_{\R^2} e^{-i (\sinh 4t)^{-1}(ay+bv)}
e^{-(\sinh 4t)^{-1}(y^2+v^2)} \,dy\,dv
= \pi (\sinh 4t) e^{-\frac{1}{4}(\sinh 4t)^{-1}(a^2+b^2)}.
$$

Combining results yields 
\begin{align*}
&\int_{\C^2} p_{2t}(z+a,w+b) e^{\frac{i}{2}(aw-bz)}
\overline{p_{2t}(z,w)} W_t(z,w)\,dz\,dw\\
&= 8^{-1} \pi^{-1} 
(\sinh 2t)^{-3}(\tanh 2t) (\sinh 4t) e^{-\frac{1}{8}(\coth 2t +
\tanh 2t)(a^2+b^2)} e^{-\frac{1}{4}(\sinh 4t)^{-1}(a^2+b^2)}.
\end{align*}
Finally using the identities $ \coth 2t+\tanh 2t = 2 \coth 4t $ and $\coth 4t
+(\sinh 4t)^{-1} = \coth 2t $ and simplifying we get
\begin{align*}
& \int_{\C^2} p_{2t}(z+a,w+b) e^{\frac{i}{2}(aw-bz)}
\overline{p_{2t}(z,w)} W_t(z,w)\,dz\,dw\\
& \hspace{1in}
=  \frac{1}{4\pi}(\sinh 2t)^{-1} e^{-\frac{1}{4}\coth 2t (a^2+b^2)}
=  p_{2t}(a,b)\, . 
\end{align*}
This proves the lemma. 
\end{proof}

\subsection{The twisted Bergman space and surjectivity of $H_t^\lambda$}
\setcounter{equation}{0}

For each $\lambda\in\R$, $\lambda\neq 0$,  we define the 
{\it $\lambda$-twisted Bergman space} by 
$$
\CB_t^\lambda(\C^{2n})= \{f\in \CO(\C^{2n}) : \|f\|_\lambda^2 =
\int_{\C^n \times \C^n} |f(\z,\w)|^2 W_t^\lambda(\z,\w)\,
 d\z\,d\w < \infty\} \, . 
$$
Clearly $\CB_t^\lambda(\C^{2n})$ is a Hilbert space of holomorphic
functions on $\C^{2n}$. 
It follows from Proposition \ref{prop:4.1} that 
$H_t^\lambda: L^2(\R^{2n})\to \CB_t^\lambda(\C^{2n})$ is an isometric 
embedding. 

\par Our goal for this subsection is to show 
that $H_t^\lambda$ is onto. We begin with 
a description of a useful orthonormal basis for 
$\im H_t^\lambda$ in terms  of  the special Hermite functions 
$\Phi_{\alpha,\beta}^\lambda (\x,\u)$ (see \cite[Section 2.3]{T}). 
For each $\alpha,\beta \in \N_0^n$, let us consider 
$$
 \wt \Phi_{\alpha,\beta}^\lambda (\z,\w) = (2\pi)^{-n} 
e^{-(2|\beta| +n)|\lambda|t}
    \Phi_{\alpha,\beta}^\lambda (\z,\w)
$$
where $\Phi_{\alpha,\beta}^\lambda (\z,\w)$ is the extension of 
$\Phi_{\alpha,\beta}^\lambda (\x,\u)$ to $\C^n \times \C^n$. The functions
$\Phi_{\alpha,\beta}^\lambda (\x,\u)$ satisfy the orthogonal relation
$$
 (\Phi_{\alpha,\beta}^\lambda *_\lambda \Phi_{\mu,\nu}^\lambda) (\x,\u) = 
   \delta_{\beta,\mu} \Phi_{\alpha,\nu}^\lambda (\x,\u)\, .  
$$

\begin{lem}
The set $\{\wt \Phi_{\alpha,\beta}^\lambda: \alpha,\beta\in \N_0^n\}$ is an 
orthonormal basis for $\im H_t^\lambda$.
\end{lem}

\begin{proof}
It is enough to prove it for $\lambda =1$ and we drop the superscript when
$\lambda=1$. As the heat kernel $p_t(\x,\u)$ is given by
$$
 p_t(\x,\u) = (2\pi)^{-n} \sum_\mu e^{-(2|\mu|+n)t} \Phi_{\mu,\mu}(\x,\u)\ ,
$$
we obtain the relation
$$
(\Phi_{\alpha,\beta} \times p_t) (\x,\u) = 
  (2\pi)^{-n} e^{-(2|\beta|+n)t} \Phi_{\alpha,\beta}(\x,\u)\, .
$$
Thus $H_t(\Phi_{\alpha,\beta}) (\z,\w) = \wt \Phi_{\alpha,\beta}(\z,\w)$
and, therefore, using Proposition \ref{prop:4.1} we obtain
\begin{align*}
& \int_{\C^{2n}} \wt \Phi_{\alpha,\beta}(\z,\w) \overline{
  \wt \Phi_{\mu,\nu}(\z,\w) } W_t(\z,\w)\, d\z \,d\w \\  
& \quad = \int_{\C^{2n}} H_t(\Phi_{\alpha,\beta})(\z,\w) \overline{
   H_t(\Phi_{\mu,\nu})(\z,\w) } W_t(\z,\w) \,d\z \,d\w \\  
& \quad = \int_{\R^{2n}} \Phi_{\alpha,\beta}(\x,\u) \overline{
   \Phi_{\mu,\nu}(\x,\u) } \,d\x \,d\u\, .   
\end{align*}
Hence $\{\wt \Phi_{\alpha,\beta}: \alpha,\beta\in \N_0^n\}$ is an 
orthonormal system in $\im H_t$. 

To show that it is an orthonormal basis for $\im H_t$, we only need to
show that 
$$
 \int_{\C^{2n}} H_t(f) (\z,\w) \overline{\wt \Phi_{\alpha,\beta}(\z,\w)}
    W_t(\z,\w) \,d\z \,d\w =0  
$$
for all $\alpha,\beta$ implies $f \equiv 0$. But the above simply means,
by Proposition \ref{prop:4.1}, that
$$
\int_{\R^{2n}} f(\x,\u) \overline{ \Phi_{\alpha,\beta}(\x,\u)} \,d\x \ d\u =0
$$
for all $\alpha,\beta$ and we know that $\{\Phi_{\alpha,\beta}: 
\alpha,\beta\in \N_0^n\}$ is an orthonormal basis for $L^2(\R^{2n})$. Hence
$f \equiv 0$ and the proof is complete.
\end{proof}

We will show that $\{\wt \Phi_{\alpha,\beta}^\lambda: \alpha,\beta\in \N_0^n\}$
is also an orthonormal basis for $\CB_t^\lambda(\C^{2n})$. Clearly this 
implies that $H_t^\lambda: L^2(\R^{2n})\to \CB_t^\lambda(\C^{2n})$ is 
onto. 

\par Note that $\wt \Phi_{\alpha,\beta}^\lambda \in \CB_t^\lambda(\C^{2n})$ 
for any $t > 0$ and 
$\{\wt \Phi_{\alpha,\beta}^\lambda: \alpha,\beta\in \N_0^n\}$ will be an 
orthonormal basis for any $\CB_t^\lambda(\C^{2n})$. 
\par As $\z = \x + i \y$ and $\w = \u + i \v$, we note that $\u \cdot \y - 
\v \cdot  \x = {\Im (\z \cdot \overline{\w})}$ is the symplectic form on 
$\R^{2n}$. Thus $\Im (\sigma \z \cdot\overline{\sigma \w}) = 
{\Im (\z \cdot\overline{\w})}$ for $\sigma\in U(n)$. 
\par We introduce the {\it twisted   
Fock space} $\CF_t^\lambda(\C^{2n})$  by 
\begin{align*}
\CF_t^\lambda(\C^{2n})=\{ &G\in  \CO(\C^{2n}): \\
&\|G\|^2=  
\int_{\C^n\times \C^n} |G(\z,\w)|^2 e^{\lambda\Im  (\z\cdot\overline{\w})} 
     e^{-\frac\lambda2 (\coth 2t\lambda) (|\z|^2+|\w|^2)} \,d\z \,d\w <
 \infty\}\, .
\end{align*}
Clearly, the prescription 
$$
U(n)\times \CF_t^\lambda(\C^{2n})\to \CF_t^\lambda(\C^{2n}), \ \ 
(\sigma, G)\mapsto 
G^\sigma; \ G^\sigma(\z,\w)= G(\sigma\z, \sigma \w)
$$ 
defines a unitary representation of $U(n)$ on $\CF_t(\C^{2n})$. 

\par The Hilbert spaces $\CB_t^\lambda(\C^{2n})$ and $\CF_t^\lambda(\C^{2n})$ 
are related through 
\begin{equation}\label{eq:b}
 F(\z,\w) \in \CB_t^\lambda(\C^{2n}) \quad \hbox{if and only if} \quad F(\z,\w)
 e^{\frac\lambda4 (\coth 2t\lambda) (\z\cdot \z+\w \cdot \w)} 
\in \CF_t^\lambda(\C^{2n}) \, .
\end{equation}

Let $T\simeq ({\S}^1)^n$ be the diagonal subgroup 
of $ U(n)$.  We write
the elements of $T$ as $ \sigma=(e^{i\varphi_1},\ldots,
e^{i\varphi_n})$.  For each $ n-$tuple of integers $ \m = (m_1,m_2,\ldots,m_n)$
let $ \chi_{\m}(\sigma) $ be the character of $T$ defined by $ \chi_{\m}(
\sigma) = e^{i \sum_{j=1}^n m_j \varphi_j }$.  For each 
$ G \in \CF_t^\lambda(\C^{2n}) $ define
$$
G_{\m}(\z,\w) = 
\int_T G(\sigma\z , \sigma \w )\overline{\chi_{\m}(\sigma)}\, d\sigma\, .
$$
As $G$ is holomorphic it is clear that $ G_{\m} = 0 $ unless $ \m $ is a
multi-index in $ \N_0^n$.  By the Fourier expansion 
$$
G(\sigma \z,\sigma \w) = \sum_{\m \in \N_0^n} G_{\m}(\z,\w)\chi_\m (\sigma)\, 
$$
and by the Plancherel theorem we have 
\begin{equation} \label{eq:4.3.2}
  \int_T |G(\sigma \z, \sigma \w)|^2 d \sigma =
       \sum_{\m \in \N_0^n} |G_{\m}(\z,\w)|^2 \ .
\end{equation}
Note that the functions $ G_{\m} $ satisfy the homogeneity condition
$$
G_{\m}(\sigma \z,\sigma \w) = \chi_{\m}(\sigma) G_{\m}(\z, \w) \ .
$$
For any $G \in \CF_t^\lambda(\C^{2n})$ we observe that, as $\Im (\z \cdot \overline{\w})
 = \Im (\sigma \z \cdot \sigma \overline{\w})$, 
\begin{align*}
& \int_{\C^{2n}} G(\z,\w) e^{\lambda \Im(\z \cdot \overline{\w} )}
   e^{-\frac{\lambda}{2} (\coth 2t \lambda)(|\z|^2+|\w|^2)} d\z\ d\w \\
& \qquad \quad  = \int_T \int_{\C^{2n}} G(\sigma \z, \sigma \w) e^{\lambda 
     \Im (\z \cdot \overline{\w})}
 e^{-\frac{\lambda}{2} (\coth 2t \lambda)(|\z|^2+|\w|^2)} d\z \ d\w \ d\sigma.
\end{align*}
In view of this and the homogeneity condition we arrive at 
the orthogonality relations  
$$
 \int_{\C^{2n}} G_\m(\z,\w) \overline{G_{\m'}(\z,\w)} 
e^{\lambda \Im(\z \cdot \overline{\w})} 
   e^{-\frac{\lambda}{2} (\coth 2t \lambda)(|\z|^2+|\w|^2)} d\z\ d\w =0\, , 
$$
whenever $\m $ and $\m'$ are different. We also note that each $G_\m$ has 
an expansion of the form 
$$ 
G_{\m}(\z,\w) = \sum_{\alpha +\beta =\m} 
c_{\alpha,\beta} \z^\alpha \w^\beta\, .
$$
Hence each $ G_{\m} $ is a polynomial.

\begin{lem} \label{lem:4.4} 
The linear span of $ P_{\alpha,\beta}(\z,\w) = \z^\alpha \w^\beta$, 
$\alpha,\beta \in \N_0^n$, is dense in $\CF_t^\lambda(\C^{2n})$ .
\end{lem} 

\begin{proof}
If $G \in \CF_t^\lambda(\C^{2n})$ is orthogonal to all $P_{\alpha,\beta}$ then
$$
 \int_{\C^{2n}} G(\z,\w) \overline{G_{\m}(\z,\w)} 
e^{\lambda \Im(\z \cdot \overline{\w})} 
   e^{-\frac{\lambda}{2} (\coth 2t \lambda)(|\z|^2+|\w|^2)} d\z\ d\w =0 
$$
for any $ \m \in \N_0^n $. In view of the homogeneity property of $G_\m$ this
means that 
$$
 \int_{\C^{2n}} |G_\m(\z,\w)|^2
e^{\lambda \Im(\z \cdot \overline{\w})} 
   e^{-\frac{\lambda}{2} (\coth 2t \lambda)(|\z|^2+|\w|^2)} d\z\ d\w =0 \ . 
$$
Hence $G_\m(\z,\w) =0$ for every $\m$ and so $G =0$ in view of 
\eqref{eq:4.3.2}. 
\end{proof}

It follows from Lemma \ref{lem:4.4} and (\ref{eq:b})
that every $F \in \CB_t^\lambda(\C^{2n})$ has the orthonormal expansion 
\begin{align} \label{eq:4.4}
  F(\z,\w) = \sum_{\m}\sum_{\alpha+\beta = \m} c_{\alpha,\beta}
   P_{\alpha,\beta}(\z,\w)
e^{-\frac\lambda4 (\coth 2t\lambda)(\z\cdot \z+\w \cdot \w)}\, .
\end{align}
The functions 
$$
\Psi_\m(\z,\w) = \sum_{\alpha+\beta = \m} 
  c_{\alpha,\beta} P_{\alpha,\beta}(\z,\w)e^{-\frac\lambda4 
(\coth 2t\lambda)(\z\cdot \z+\w \cdot \w)}                             
$$
are orthogonal in $ \CB_t^\lambda(\C^{2n}) $ but not orthogonal in any other 
$\CB_s^\lambda( \C^{2n})$ when $s \ne t$. Another crucial property of these
functions is proved in the next lemma.

\begin{lem} \label{lem:4.5}  
All the functions $  \Psi^{\m}_{\alpha,\beta}(\z,\w) = 
 P_{\alpha,\beta}(\z,\w)e^{-\frac\lambda4 
(\coth 2t\lambda)(\z\cdot \z+\w \cdot \w)}$ belong to the image 
$ \im H_t^\lambda $ of the heat kernel transform.
\end{lem} 

\begin{proof} We may restrict ourselves to the case of $\lambda=1$. 
It will suffice to show that for each pair $ \alpha, \beta \in \N_0^n$ 
there exists
a function $f_{\alpha,\beta} \in L^2(\R^{2n})$ such that 
$$
H_t(f_{\alpha,\beta})(\z,\w)= (f_{\alpha,\beta}\times p_t)^\sim (\z,\w)
 = \z^\alpha \w^\beta e^{-\frac14 (\coth 2t)
 (\z^2+\w ^2)}\, .   
$$
As both sides are holomorphic it is enough to prove this for $\z = \x$ and 
$\w = \u$ where $ \x,\u \in \R^n.$ Thus we need to solve the equation
\begin{equation}\label{eq:4.5}
(f_{\alpha,\beta}\times p_t)(\x,\u) = \x^\alpha \u^\beta p_{2t}(\x,\u)\, .
\end{equation}
In the sequel it will be convenient to identify 
$\R^{2n}$ with $\C^n$ via $ z = \x+i\u$. Then 
$\x^\alpha \u^\beta = 2^{-|\alpha|}(2i)^{-|\beta|}
(z+\overline{z})^\alpha (z-\overline{z})^\beta $. It is then sufficient to
solve the equation
$$
(f_{\alpha,\beta}\times p_t)(z) = z^\alpha \overline{z}^\beta p_{2t}(z)
$$
where $ p_t(z) = p_t(\x,\u).$ We solve this equation using properties of the
Weyl transform.

Recall that the Weyl transform ${\W}(f)$ of a function  $f \in L^1(\C^n)$,
is defined to be the bounded operator on $L^2(\R^n)$ given by
$$
{ \W}(f)\varphi(\xi) = \int_{\C^n} f(z) \pi(z) \varphi(\xi)\, 
dz\qquad (\xi\in \R^n) 
$$   
where $\pi(z) = \pi_1(z,0)$ and $\pi_1$ is the Schr\"odinger 
representation of the Heisenberg group $\H$ with parameter $\lambda =1$
(see \cite[Section 2.2]{T}). Then for $f \in L^1 \bigcap L^2(\C^n)$, 
${\W}(f)$ 
is a Hilbert-Schmidt operator and ${\W}$ extends to $L^2(\C^n)$ as an 
isometry
onto the space of Hilbert-Schmidt operators. Moreover ${\W}(f \times g) = 
{\W}(f){\W}(g)$ and ${\W}(p_t) = e^{-t H}$. Here $H$ denotes 
the Hermite operator
$$
  H = (-\Delta + |\xi|^2) = \frac{1}{2}\sum_{j=1}^n(A_j A_j^*+ A_j^* A_j), 
$$ 
in which $A_j = - \frac{\partial}{\partial \xi_j} +\xi_j$ and $A_j^* = 
 \frac{\partial}{\partial \xi_j} + \xi_j$ are the creation and annihilation
operators. The eigenfunctions of $H$ are the Hermite functions $\Phi_\alpha$.
They satisfy
$$
 A_j \Phi_\alpha = (2 \alpha_j+2)^{\frac12} \Phi_{\alpha + e_j},  \qquad 
 A_j^* \Phi_\alpha = (2 \alpha_j)^{\frac12} \Phi_{\alpha - e_j}
$$
where $e_j$ are the coordinate vectors. 
Given a bounded linear operator $T$ on $L^2(\R^n)$, define the derivations
$$
 \delta_j T = [A_j^*,T] = A_j^* T - T A_j^*, \qquad   
 \overline{\delta}_j T = [T, A_j] = T A_j - A_j T\, .  
$$ 
Then it can be shown that (see \cite{T2})
$$
  {\W}(z_j f) = \delta_j {\W}(f), \quad \hbox{and} \quad 
  {\W}(\overline{z}_j f) = \overline{\delta}_j {\W}(f)\, .
$$ 
By iteration we obtain
$$
{\W}(z^\alpha \overline{z}^\beta f) = \delta^\alpha \overline{\delta}^\beta
   {\W}(f)
$$
where $\delta^\alpha \overline{\delta}^\beta$ are defined in an obvious way.

Returning to our equation \eqref{eq:4.5}, we take the Weyl transform on
both sides and obtain that 
$$
  {\W}(f_{\alpha,\beta}) e^{-t H} = \delta^\alpha \overline{\delta}^\beta 
e^{-2 t H} \, .
$$
Testing against the Hermite basis it is easy to see that the densely 
defined operator
$$
  T = (\delta^\alpha \overline{\delta}^\beta e^{-2t H}) e^{t H}
$$
extends to the whole $L^2(\R^n)$ as a Hilbert-Schmidt operator. Hence,
$T = {\W}(f_{\alpha,\beta})$ for some $f_{\alpha, \beta} \in L^2(\C^n)$. 
This completes the proof of the lemma. 
\end{proof}  

\begin{thm} \label{thm:4.6} Let $t>0$ and $\lambda\in \R$, $\lambda\neq 0$. 
Then the $\lambda$-twisted heat kernel transform 
$H_t^\lambda: L^2(\R^{2n})\to \CB_t^\lambda(\C^{2n})$ is an 
isometric isomorphism. Moreover, $\{\wt \Phi_{\alpha,\beta}^\lambda:  
\alpha,\beta\in \N_0^n\}$ is an orthonormal basis for $\CB_t^\lambda(\C^{2n})$.
\end{thm} 

\begin{proof} 
As usual we restrict our attention to the case $\lambda =1$. All what is 
left to show is that $H_t$ is onto. Suppose that $F\in \CB_t(\C^{2n})$ 
is orthogonal to all $\wt \Phi_{\alpha,\beta}$. 
We have to verify that $F\equiv 0$. The function 
$$
   G(\z,\w) = F(\z,\w) e^{\frac14 (\coth 2t)(\z\cdot \z+\w \cdot \w) }
$$
is orthogonal in $\CF_t$ to all functions of the form
$$
  f \times p_t(\z,\w) e^{-\frac14 (\coth 2t)(\z\cdot \z+\w \cdot \w) }.
$$
In view of Lemma \ref{lem:4.5}, $G$ is orthogonal to all $P_{\alpha,\beta}$.
Hence by Lemma \ref{lem:4.4} we get $G= 0$ and so $F =0$ as desired.
\end{proof}

We conclude this subsection with a proof of the uniqueness of the 
weight function $W_t^\lambda$. 

\begin{lem}\label{lem=unique} $W_t^\lambda$ is the unique non-negative measurable 
weight function for the $\lambda$-twisted Bergman space
$\CB_t^\lambda(\C^{2n})$.  
\end{lem}

\begin{proof} In view of (\ref{eq:b}),  the statement is 
equivalent to the assertion that
\begin{equation}\label{bird} \CW_t^\lambda(\w,\z)=
e^{\lambda \Im(\z \cdot \overline{\w})} 
   e^{-\frac{\lambda}{2} (\coth 2t \lambda)(|\z|^2+|\w|^2)}
\end{equation}
is the unique weight function for the twisted Fock space
$\CF_t^\lambda(\C^{2n})$. This will be verified in the sequel. 
\par We may restrict ourselves to the notationally 
convenient case $n=1$, $\lambda=1$ and drop all
sub and superscripts involving $\lambda$. 
Let $\CU_t:\C^2\to \R_{\geq 0}$ be a measurable function such that 
\begin{equation} \label{tra}
\int_{\C^2} f(z,w)\oline {g(z,w)}\  \CW_t(z,w) \, dz\,  dw
=\int_{\C^2} f(z,w)\oline {g(z,w)}\ \CU_t(z,w) \, dz\,  dw
\end{equation}
holds for all $f,g\in \CF_t(\C^2)$. We have to show that 
$\CW_t=\CU_t$ almost everywhere. Recall from 
Lemma \ref{lem:4.5}  that 
all polynomials $z^m w^n$ lie in $\CF_t(\C^2)$. In particular
the constant function belongs to  $\CF_t(\C^2)$
and (\ref{tra}) implies that $\CU_t$ is integrable.  
\par Let us introduce polar coordinates 
on $\C^2$ by $(z,w)=(re^{i\phi}, se^{i\theta})$.  
Consider the Fourier expansions of $\CW_t$ and 
$\CU_t$ given by 
$$\CW_t(re^{i\phi}, se^{i\theta})=\sum_{m,n\in \Z} 
a_{m,n}(r,s) e^{im\phi} e^{in\theta}$$
and 
$$\CU_t(re^{i\phi}, se^{i\theta})=\sum_{m,n\in \Z} 
b_{m,n}(r,s) e^{im\phi} e^{in\theta}\, . $$
Identity (\ref{tra}) applied to $f=g=z^kw^l$
yields the estimates
\begin{align}\label{tro}
& \int_0^\infty \int_0^\infty r^{2k+1} s^{2l+1} |a_{m,n}(r,s)| \, dr\,  ds 
\leq \| z^k w^l\|^2_{\CF_t(\C^2)}\\
\label{tro1}& \int_0^\infty \int_0^\infty r^{2k+1} s^{2l+1} |b_{m,n}(r,s)| \, dr\,  ds 
\leq \| z^k w^l\|^2_{\CF_t(\C^2)}
\end{align}
for all $m,n\in \Z$.

\par We finish the proof and show $c_{m,n}=a_{m,n}-b_{m,n}=0$ for all 
$m,n\in \Z$. 
In fact for $f=z^{m_1} w^{n_1}$ and $g=z^{m_2}w^{n_2}$ for $m_1, m_2, n_1,n_2\in \N_0$ 
we obtain from (\ref{tra}) that 
\begin{equation}\label{tri}\int_0^\infty \int_0^\infty 
 r^{m_1 +m_2 +1} s^{n_1+n_2+1}  c_{m_2-m_1, n_2-n_1}(r,s) \, dr\, ds=0\, .
\end{equation} 
Note that the integral on the left is absolutely convergent by 
(\ref{tro})-(\ref{tro1}). 
Fix now $m,n\in \Z$. Reformulating (\ref{tri}) reads 

\begin{equation} \label{whale}\int_0^\infty \int_0^\infty 
 r^{|m| + 2k+1} s^{|n|+2l+1}  c_{m, n}(r,s)\, dr\, ds
= 0 \end{equation}
for all $k,l\in \N_0$.
In view of (\ref{tro})-(\ref{tro1}), we have the estimate 
\begin{equation} \label{shark}\int_0^\infty \int_0^\infty 
 r^{|m| + 2k+1} s^{|n|+2l+1}  |c_{m, n}(r,s)|\, dr\, ds
\leq 2 \|z^{|m|+k} w^{|n| +l}\|^2_{\CF_t(\C^2)}+C
\end{equation}
with $C=\int_{|z|<1, |w|<1}
 (\CW_t(z,w)+ \CU_t(z,w))\, dz\, dw>0$ a constant 
independent of $m,n$.
  
\par Denote by $\CR_+=\{ \zeta\in \C: 
 \Re \zeta>0\}$ the right halfplane. Let us 
recall the elementary fact that a bounded 
holomorphic function $f: \CR_+\to \C$ 
which vanishes on $\alpha +\beta \N_0$ for some 
$\alpha\geq 0$, $\beta>0$ is identically zero 
(see \cite{K}, Lemma A.1  for a proof). 

\par  The explicite formula for $\CW_t$ in 
(\ref{bird}) yields a crude but sufficient estimate
for the norm of monomials: there exists  constants 
$c,\gamma>0$ such that for all $k,l\in \N_0$
one has 
\begin{equation}\label{starbird}
\|z^k w^l\|^2 \leq c\cdot e^{\gamma(k+l)}\, .
\end{equation} 
\par 
Now define the function 
\begin{align*} F_{m,n}:&\CR_+\times \CR_+\to \C,\\
&(\zeta_1,\zeta_2)\mapsto  
e^{-3\gamma(\zeta_1+\zeta_2)}\int_0^\infty \int_0^\infty 
 r^{|m| +2\zeta_1+1 } s^{|n|+2\zeta_2+1}  c_{m, n}(r,s)\, dr\, ds\, .
\end{align*}
It is a consequence of (\ref{shark}) and  (\ref{starbird}) that 
$F_{m,n}$ is bounded and holomorphic on $\CR_+\times \CR_+$. 
As $F_{m,n}|_{\N\times \N}=0$ by (\ref{whale}), we conclude that 
$F_{m,n}= 0$. But then $c_{m,n}=0$ by the properties of the 
Mellin transform. 
\end{proof}

\subsection{The inversion formula for $H_t^\lambda$}
\setcounter{equation}{0}

We conclude this section by proving a formula for the 
inverse map of the $\lambda$-twisted heat kernel transform 
$H_t^\lambda: L^2(\R^{2n})\to \CB_t^\lambda (\C^{2n})$.  
It is in the nature of the problem that $(H_t^\lambda)^{-1}$ can only 
be defined nicely on a dense subspace of $\CB_t^\lambda(\C^{2n}).$ 
The precise statement is as follows:

\begin{thm} \label{thm:4.7} The inverse of $H_t^\lambda: L^2(\R^{2n})\to 
\CB_t^\lambda(\C^{2n})$ is given by 
$$
(H_t^\lambda)^{-1}(F)=
\lim_{s\to 0^+} F_s\qquad (F\in \CB_t^\lambda (\C^{2n}))\, ,
$$  
where 
$$
F_s(\a,\b) = \int_{\C^{2n}} F(\z+\a,\w+\b)
 e^{\frac{i\lambda}{2}(\a\cdot\w-\b\cdot\z)}
\overline{p_{t+s}^\lambda(\z,\w)} W_t^\lambda(\z,\w) \,d\z\, d\w \, . 
$$
\end{thm}

\begin{proof} As before we only need to handle the case 
of $\lambda=1$. 

\par Let $F\in \CB_t(\C^{2n})$. Since the space $ \CB_t(\C^{2n})$ is 
twisted-translation invariant, it is clear that
the function
$$ 
(\tau(-\a,-\b)F)(\z,\w)=
F(\z+\a,\w+\b)e^{\frac{i}{2}(\a\cdot\w-\b\cdot\z)}
$$
belongs to $ \CB_t(\C^{2n}).$ Hence, by Cauchy-Schwarz ineqaulity, the 
integral defining $ F_s $
converges. According to Theorem \ref{thm:4.6} we have 
$ F = H_t(f)=(f\times p_t)^\sim $ for some  $f \in L^2(\R^{2n}) $. 
It is easy
to see that $ F_s \in L^2(\R^{2n})$ and that $F_s$ converges to $f$.  In fact,
we have 
\begin{align*}
F_s(\a,\b) & = 
\int_{\C^{2n}} (\tau(-\a,-\b)H_t(f))(\z,\w) \overline{H_t(p_s)(\z,\w)} 
W_t(\z,\w) \,d\z \,d\w \\
& = \int_{\C^{2n}} H_t(\tau(-\a,-\b)f)(\z,\w) \overline{H_t(p_s)(\z,\w)} 
W_t(\z,\w) \,d\z \,d\w \\
& = \int_{\R^{2n}}(\tau(-\a,-\b)f)(\x,\u) p_s(\x,\u) \,d\x \,d\u\, .
\end{align*}
As $(p_s)_{s>0}$ is a Dirac sequence, it therefore follows that 
$$
F_s(\a,\b)\to (\tau(-\a,-\b)f)(0,0)=f(\a,\b)
$$
for $s\to 0^+$. 
This proves the theorem. 
\end{proof}

\renewcommand{\theequation}{\thesection.\arabic{equation}}

\section{The image of $\CH_t$ as a direct integral}
\setcounter{equation}{0}

The goal of this section is to give a natural $\H$-equivariant identification 
of the image of the 
heat kernel transform $\CH_t: L^2(\H)\to \CO(\H_\C)$ with 
a direct integral of twisted Bergman-spaces. 

\par We set $\R^\times= \R\bs\{0\}$. For each $\lambda\in \R^\times$ we write 
$\la \cdot, \cdot\ra_\lambda$ for the inner product on 
$\CB_t^\lambda(\C^{2n})$. 
Recall the orthonormal basis $\{ \tilde \Phi_{\alpha, \beta}^\lambda: \alpha, 
\beta\in \N_0^n\}$
of $\CB_t^\lambda(\C^{2n})$ from Theorem 4.6. 

\par We now introduce a measurable structure on 
$\coprod_{\lambda\in \R^\times} \CB_t^\lambda(\C^{2n})$. 
By a {\it section} $s$ of  $\coprod_{\lambda\in \R^\times} 
\CB_t^\lambda(\C^{2n})$ we understand an assignment 
$$
s: \R^\times\to \coprod_{\lambda\in \R^\times} \CB_t^\lambda(\C^{2n}), \ \ 
\lambda\mapsto s_\lambda\in \CB_t^\lambda(\C^{2n})\, .
$$
We declare a section $s=(s_\lambda)$ to be {\it measurable} if for all 
$\alpha, \beta\in \N_0^n$ the map
$$
\R^\times\to \C, \ \ \lambda\mapsto \la s_\lambda, 
  \tilde \Phi_{\alpha, \beta}^\lambda\ra_\lambda
$$
is measurable. With that we can define a direct integral of Hilbert spaces 
by 
\begin{align*}
\int_{\R^\times}^\oplus \CB_t^\lambda(\C^{2n}) \, e^{2 t\lambda^2}d\lambda=
\{ s: \R^\times \to \coprod_{\lambda\in \R^\times} &\CB_t^\lambda(\C^{2n}):  
s\ \text{measurable}\, ,\\
&\|s\|^2=\int_{\R^\times} \|s_\lambda\|_\lambda^2  
\,e^{2t\lambda^2}d\lambda<\infty \}\, . 
\end{align*}
Recall the unitary representation $\tau^\lambda$ of $\H$ on
$\CB_t^\lambda(\C^{2n})$
from Subsection 4.1. We then obtain a unitary representation 
$\int_{\R^\times}\tau^\lambda\, d\lambda$ on  
$\int_{\R^\times}^\oplus \CB_t^\lambda(\C^{2n}) \, e^{2 t\lambda^2}d\lambda$
by 
$$
\left(\int_{\R^\times}\tau^\lambda\, d\lambda\right)(h)(s)
=(\tau^\lambda(h)s_\lambda)_\lambda
$$
for $h\in \H$ and $s=(s_\lambda)$ a square integrable section. 

\par In our next step we will identify $\im \CH_t$ with our direct integral 
from above. For that let $f\in S(\H)$ be a Schwartz function. Then 
$\CH_t(f)=(k_t*f)^\sim$ and from $(f*k_t)^\lambda = e^{-t\lambda^2} 
f^\lambda*_\lambda p_t^\lambda$ it hence follows that 
\begin{equation}\label{drei}
(\CH_t(f))^\lambda =e^{-t\lambda^2} H_t^\lambda(f^\lambda)\, .
\end{equation}

\begin{thm} 
Let $t>0$. The map 
$$
\CJ_t: S(\H)\to \int_{\R^\times}^\oplus \CB_t^\lambda(\C^{2n}) \, 
e^{2 t\lambda^2}d\lambda, \ \ f\mapsto \left((\CH_t(f))^\lambda\right)_\lambda
$$
extends to an $\H$-equivariant unitary equivalence 
$$(\tau, L^2(\H))\simeq \left(\int_{\R^\times}\tau^\lambda\,d\lambda\,,\, 
\int_{\R^\times}^\oplus \CB_t^\lambda(\C^{2n}) \, e^{2 t\lambda^2}d\lambda
\right)\,.$$
\end{thm}
\begin{proof} Let $f\in S(\H)$. Then 
$$\|f\|^2=\int_{\R^{2n+1}} |f(\x,\u,\xi)|^2\, d\x\, d\u\, d\xi
=\int_{\R^{2n}} \int_\R |f^\lambda(\x,\u)|^2\, d\x \, d\u\,d\lambda .$$ 
By Theorem \ref{thm:4.6} we have for each $\lambda$ that 
$$\int_{\R^{2n}} |f^\lambda(\x,\u)|^2 \, d\x\,d\u  =
\|H_t^\lambda(f^\lambda)\|^2\ .$$
Thus it follows from (\ref{drei}) that $\CJ_t$ extends to an isometric 
embedding 
$$
\CJ_t: L^2(\H)\to \int_{\R^\times}^\oplus \CB_t^\lambda(\C^{2n}) \, 
e^{2 t\lambda^2}d\lambda\, ,
$$
denoted by the same symbol. The discussion leading up to (\ref{eq:a}) 
shows that $\CJ_t$ is $\H$-equivariant. 
\par It remains to show that $\CJ_t$ is onto. For that observe if
$$
f(\x,\u,\xi)=F(\x,\u)\varphi(\xi)
$$
for Schwartz functions $F\in S(\R^{2n})$, $\varphi\in \S(\R)$, then 
$$
(\CH_t(f))^\lambda(\x,\u)=\hat\varphi(\lambda) e^{-t\lambda^2} 
H_t^\lambda(F)(\x,\u)\, .
$$
From that the surjectivity of $\CJ_t$ easily follows. 
\end{proof}

\renewcommand{\theequation}{\thesubsection.\arabic{equation}}

\section{The image of $\CH_t$ as a sum of weighted Bergman spaces}
\setcounter{equation}{0}

In this section we prove the main result of this paper: 
$\im \CH_t=\CB_t^+(\H_\C)\oplus \CB_t^-(\H_\C)$ is a direct sum of two 
weighted 
Bergman spaces. Very surprisingly, the corresponding weight functions 
$W_t^+$ and $W_t^-$ attain also negative values
(see the phenomenon explained in Example 3.1). 

\par We will begin our discussion by showing that 
$\im \CH_t$ is not a weighted Bergman space corresponding to a non-negative
weight function. This will lead naturally to the definition of the 
partial weight functions $W_t^+$ and $W_t^-$ and to a proof of the 
main theorem.  

\subsection{Non-existence of a non-negative weight function} 
\setcounter{equation}{0}

The goal of this subsection is to 
discuss the non-existence 
of a  non-negative weight function $W_t$ on 
$\H_{\C}$ such that 
\begin{equation}\label{vier}
\|f\|^2= \int_{\H_{\C}} |\CH_t(f)(z)|^2  W_t(z) \, dz
\end{equation}
holds for all $f\in  L^2(\H)$. In other words, $\im \CH_t$ is not a weighted 
Bergman space corresponding to a non-negative weight function $W_t$.
Subject to the natural assumption that 
$W_t$ is $\H$-invariant, this will be established in Theorem  \ref{thm:5.3}
below. 

\par Recall that we identify $\H_{\C}$ with $\C^n \times \C^n \times \C$.  
If $\z =\x+i\y, \w = \u+i\v , \zeta = \xi+i\eta $ then $(\z,\w,\zeta) =
he^{iX}$ with $X = (\y,\v,\eta+\frac{1}{2}(\x\cdot\v-\u\cdot\y))$ and 
$h =(\x,\u,\xi)$.  
\par Suppose that (\ref{vier}) holds. As $\CH_t$ is $\H$-equivariant, it
is natural to assume that $W_t(he^{iX}) = W_t(e^{iX})$ for all $h \in \H$. 
In coordinates $(\z,\w,\zeta) $ this means that
\begin{align} \label{eq:5.1}
W_t(\x+i\y,\u+i\v,\xi+i\eta) = W_t\left(i\y,i\v,i\eta+
       \frac{i}{2}(\x\cdot \v-\u\cdot \y)
\right).
\end{align} 
Thus the weight function is uniquely determined by its restriction to
$(i\y,i\v,i\eta)$. Furthermore  $W_t $ is independent of the $\xi$
variable. Hence (\ref{vier}) reads as 
\begin{equation} \label{eq:5.2}
\|f\|^2= \int_{\H_{\C}} |\CH_t(f)(\z,\w,\zeta)|^2  W_t(\z,\w,i\eta) 
\,d\z\,d\w\,d\zeta 
\end{equation} 

\begin{prop} \label{prop:5.1}
Let $W_t(\z,\w,i\eta)$ be a non-negative measurable  function on $\H_\C$. If 
\eqref{eq:5.2} holds for all $f\in L^2(\H)$, then it is necessary 
that $W_t$ satisfies 
\begin{equation}\label{eq:5.3}  
W_t^\lambda(\z,\w)  = e^{- 2t \lambda^2}
  \int_\R  e^{2\lambda \eta}
    W_t(\z,\w,i\eta)  \,d\eta
\end{equation}
for all $\lambda\in \R^\times$ and 
$W_t^\lambda$ the function given in \eqref{eq:4.3}.
\end{prop}

\begin{proof} Write 
$$
{\mathcal W}_t^\lambda(\z,\w)=e^{- 2t \lambda^2}
  \int_\R  e^{2\lambda \eta}
    W_t(\z,\w,i\eta) \,d\eta\ .
$$
We have to show that $W_t^\lambda={\mathcal W}_t^\lambda$. 
\par It follows from (2.2.2) that  
$$
\int_{\R} k_t^\sim(\z, \w, \xi+i\eta) e^{i\lambda \xi}\, d\xi 
= e^{\lambda \eta} e^{-t\lambda^2} p_t^\lambda(\z, \w)\ .
$$
An easy calculation shows that
\begin{align}\label{eq:5.4}
\int_{\R} \CH_t(f) (\z,\w, \xi+i\eta) e^{i\lambda \xi}\, d\xi 
& = e^{\lambda \eta} e^{-t\lambda^2} (f^\lambda *_\lambda
  p_t^\lambda)(\z,\w) \\
& =e^{\lambda \eta} e^{-t\lambda^2} H_t^\lambda(f^\lambda)(\z,\w)\ . \notag
\end{align}
Therefore, upon applying Plancherel theorem in the $\xi$-variable,
the equation \eqref{eq:5.2} becomes
$$ 
\|f\|^2=\int_\R \int_{\C^{2n}}\int_\R  |H_t^\lambda(f^\lambda)(\z,\w)|^2  
e^{-2t\lambda^2} e^{2\lambda\eta} W_t(\z,\w, i\eta)
\,d\eta\, d\x\,d\u\,d\y\,d\v\, d\lambda\ .
$$
Here we applied Fubini's theorem which is justified 
as $W_t$ is by assumption non-negative. Employing the definition 
of ${\mathcal W}_t$ we therefore get 
$$ 
\int_\R \int_{\R^{2n}}|f^\lambda(\x,\u)|^2 \,d\x \,d\u\,d\lambda
=\int_{\R} \int_{\C^{2n}} |H_t^\lambda(f^\lambda)(\z,\w)|^2 
{\mathcal W}_t^\lambda (\z,\w)
\,d\x\,d\u\,d\y\,d\v\, d\lambda \, .
$$ 

Let now $ \varphi $ be a Schwartz class function on $\R$ with unit $L^2$-norm
and define $ f $ by $ f(\x,\u,\xi) = \wh\varphi(\xi)F(\x,\u)$
with $F\in L^2(\R^{2n})$.  Then $f^\lambda(\x, \u)=\varphi(\lambda) F(\x,\u)$
and $H_t^\lambda(f^\lambda)=\varphi(\lambda) H_t^\lambda(F)$.
For such $f$ the above displayed equation becomes
\begin{equation}\label{sechs}
\int_{\R^{2n}}|F(\x,\u)|^2 \,d\x \,d\u
=\int_{\R} \int_{\C^{2n}} |\varphi(\lambda)|^2
|H_t^\lambda(F)(\z,\w)|^2 {\mathcal W}_t^\lambda (\z,\w)
\,d\x\,d\u\,d\y\,d\v\, d\lambda \, .
\end{equation} 
From (\ref{sechs}) it is easy to see that for every $\lambda\neq 0$ and 
all  $F\in L^2(\R^{2n})$
$$  \int_{\R^{2n}} |F(\x,\u)|^2  \,d\x\,d\u=
\int_{\C^{2n}} |H_t^\lambda(F)(\z,\w)|^2 
{\mathcal W}_t^\lambda(\z, \w)  \,d\x \,d\u\,d\y\,d\v \, . 
$$
By Lemma  \ref{lem=unique}, the weight function
$ \CW_t^\lambda $ is given by \eqref{eq:4.3}.
\end{proof}

\begin{thm} \label{thm:5.3} 
There is no non-negative left $\H$-invariant weight function $W_t$ for which $\eqref{eq:5.2}$ 
holds for all $f \in L^2(\H)$, i.e. $\im \CH_t$ is not a weighted 
Bergman spaces corresponding to a left $\H$-invariant non-negative weight function.  
\end{thm} 

\begin{proof}
By \eqref{eq:5.1}, $W_t$ is uniquely determined by its restriction to 
$(i\y,i\v,i\eta)$. By \eqref{eq:4.3} and \eqref{eq:5.3},
$$  
\int_\R  e^{2\lambda\eta} W_t(i\y,i\v,i\eta) \,d\eta =
       e^{2 t \lambda^2} p_{2t}^\lambda (2\y, 2\v)\qquad (\lambda\in 
\R^\times)\, .$$
If $W_t$ were non-negative, then for fixed $\y,\v$ and $\lambda $ the function 
$\eta\mapsto e^{2\lambda\eta} W_t(i\y,i\v,i\eta)$ would belong to 
$L^1(\R)$. Consequently, we would have 
\begin{equation}\label{1000}
  \int_\R e^{2(\lambda+is)\eta} W_t(i\y,i\v,i\eta) \,d\eta =
       e^{2 t (\lambda+i s)^2} p_{2t}^{\lambda+is} (2\y, 2\v)\, .
\end{equation}
The left hand side of (\ref{1000}) would be holomorphic in $\lambda+is$ since 
for every $n\in \N_0$ there exists an $\eps>0$ such that 
$|\eta|^n e^{2\lambda\eta} W_t(i\y,i\v,i\eta) \le 
 e^{2\lambda\eta+\varepsilon |\eta|} W_t(i\y,i\v,i\eta)$. 
However, the right side of (\ref{1000}) is holomorphic only for 
$\lambda \ne 0$. If $\lambda =0$, it becomes
$$
p_{2t}^{is} (2 \y,2\v) = c_n\left(\frac s{\sin (2 s t)} \right)^n
    e^{-s (\cot 2 s t) (\y^2+\v^2)},
$$
which has an essential singularity at the points $s \in \Z^\times (\pi/t)$. 
Therefore there is no non-negative $W_t$ that will satisfy
\eqref{eq:5.3} or \eqref{eq:5.2}.
\end{proof}  

\subsection{The partial weight functions $W_t^+$ and $W_t^-$} 
\setcounter{equation}{0}
Recall the twisted weight function $W_t^\lambda$ from (\ref{eq:4.3}). 
\par Let $\lambda>0$ and define a function $W_t^+$ on $\H_\C$ by 
\begin{equation}\label{1001}
W_t^+(\z,\w,\zeta) = \int_\R
e^{2t(\lambda+\frac{i}{2}s)^2} e^{-2\eta(\lambda +\frac{i}{2}s)} 
W_t^{\lambda+\frac{i}{2}s}(\z,\w) \,ds . 
\end{equation}
It is easy to see that $W_t^+$ is well-defined. Notice that 
$W_t^+$ does not depend on $\xi$. 
In Proposition \ref{prop:5.3} below 
we will show that $W_t^+$ is  independent of the choice of $\lambda>0$. 

\begin{prop} \label{prop:5.3} The function $W_t^+$ satisfies 
the following properties: 
\item{{\rm (i)}} $W_t^+$ is independent of the choice
of $\lambda>0$. In particular, 
$$W_t^+(\z,\w,\zeta) = \lim_{\lambda\to 0^+} \int_\R 
e^{2t(\lambda+\frac{i}{2}s)^2} e^{-2\eta(\lambda+\frac{i}{2}s)} 
W_t^{\lambda+\frac{i}{2}s}(\z,\w) \,ds .$$
\item{{\rm (ii)}} Let $a>0$ and $Q\subseteq \C^{2n}$ be a compact set. 
Then there exists a constant $C=C(Q,a)>0$ such that 
for all $\eps\in [a^{-1}, a]$ and $\xi\in \R$ 
$$\sup_{(\z,\w)\in Q} \int_\R \left|e^{2\eps \eta} W_t^+(\z,\w,\xi+i\eta)\right|\, d\eta\leq C\, .$$
\item{{\rm (iii)}} $W_t^+$ satisfies \eqref{eq:5.3} with $\lambda>0$, i.e.
\begin{equation}\label{1002}
W_t^\lambda (\z,\w) =e^{-2t\lambda^2} \int_\R  e^{2\eta \lambda} W_t^+(\z,\w,i\eta) \,d\eta
\end{equation}
for $\lambda>0$.
\item{{\rm (iv)}} $W_t^+$ is real valued and left $\H$-invariant. 
\end{prop}

\begin{proof} (i) Let $\lambda>0$. We have to show that 
$$W_t^+(\z,\w,\zeta) = \int_\R 
e^{2t(\lambda+\frac{i}{2}s)^2} e^{-2\eta(\lambda +\frac{i}{2}s)} 
W_t^{\lambda+\frac{i}{2}s}(\z,\w) \,ds $$
is independent of the choice of $\lambda>0$. 
This will be a consequence of  Cauchy's theorem. 
Indeed, let us denote the right hand side by  $I(\lambda)$.
For $R > 0$ and $\lambda_2 > \lambda_1 >0$, let $\Gamma_R$ be the contour
consisting four lines, $\Gamma_R(\lambda_1):=\{\lambda_1+is/2: -2R < s <2R\}$,
$\gamma_{-R}=\{\lambda - iR: \lambda_1 \le \lambda \le\lambda_2\}$,
$\Gamma_R(\lambda_2) = \{\lambda_2+is/2: -2R < s <2R\}$ and $\gamma_{R} =
\{\lambda + iR: \lambda_1 \le \lambda \le\lambda_2\}$, going counterclockwise.
As $R \to \infty$, the integral on $\Gamma_R(\lambda)$ becomes $I(\lambda)$.
Cauchy's theorem shows that 
$$
\int_{\Gamma_R} e^{-2\eta z}e^{2tz^2} W_t^z(\z,\w) \,dz = 0.
$$
It is easy to see that $|\sinh (\lambda+iR)t| \ge \sinh (\lambda t)$
and $|\cosh (\lambda+iR)t| \le \cosh (\lambda t)$. Thus,
$$
|p_{2t}^{\lambda\pm iR}(2\y,2\v)| \le
     \left(\frac{\lambda+R}{\sinh \lambda t} \right)^n
       e^{(\lambda+R)\coth(\lambda t) (|\y|^2+|\v|^2)}.
$$
Together with $|e^{2t(\lambda \pm iR)^2}| = e^{2t \lambda^2} e^{-2tR^2}$, this
shows that the integrals on $\gamma_{-R}$ and on $\gamma_{R}$ go to zero
as $R\to +\infty$. Thus, taking $R\to \infty$ shows that $I(\lambda_1) =
I(\lambda_2)$. This completes the proof of (i). 

\par (ii) It follows from (i) that 
$W_t^+$ satisfies the bound
$$
|W_t^+(\z,\w, \xi+i\eta)| \le e^{- 2\eta \lambda} e^{2t \lambda^2}
    \int_\R e^{- \frac{1}{2} t s^2}
       \left|W_t^{\lambda+\frac{i}{2}s}(\z, \w)\right |\,ds\ .
$$
for any $\lambda>0$. Notice that the integral on the right is 
independent of $\eta$.  Thus if we let $\lambda >\epsilon$ if $\eta > 0$
and $\lambda <\epsilon$ if $\eta < 0$, we see that
$\eta\mapsto e^{2 \varepsilon \eta} W_t^+(\z,\w,\xi+i\eta)$
is integrable.
This implies (ii). 

\par (iii) This is immediate from the definition (\ref{1001}) and Fourier 
inversion (which is justified by (ii)).  In fact, we have 
$$
W_t^+(\z,\w,\zeta)=e^{-2\eta\lambda}\int_\R e^{-i\eta s} 
e^{2t(\lambda+\frac{i}{2}s)^2} W_t^{\lambda+\frac{i}{2}s} (\z, \w) \, ds
$$
and so 
$$
\int_\R e^{2\lambda\eta}W_t^+(\z,\w,\xi+i\eta)e^{i\eta s}\, d\eta=
e^{2t (\lambda+\frac{i}{2}s)^2} W_t^{\lambda+\frac{i}{2}s}(\z, \w)\ .
$$
Setting $s =0$ gives the the stated result.

\par (iv) We first show that $W_t^+$ is real valued. 
In fact, taking the conjugate of the integral (\ref{1001})  and then changing variable 
$s \to -s$ shows that the weight function $W_t^+$ is real. 
Finally, the fact that $W_t^\lambda$ is twisted-translation invariant 
forces that $W_t^+$ is left $\H$-invariant. 
\end{proof}
The function $W_t^+$ has a natural counterpart $W_t^-$. 
For $\lambda<0$ we define $W_t^-$ by 
\begin{equation}\label{1003}
W_t^-(\z,\w,\zeta) = \int_\R
e^{2t(\lambda+\frac{i}{2}s)^2} e^{-2\eta(\lambda +\frac{i}{2}s)} 
W_t^{\lambda+\frac{i}{2}s}(\z,\w) \,ds . 
\end{equation}

It is more or less obvious that $W_t^-$ satisfies the 
same properties as $W_t^+$ listed in Proposition (\ref{prop:5.3}), i.e.
$W_t^-$ is independent of the choice of $\lambda<0$ etc. In fact, a 
simple change of variable in the integral and the fact that 
$p_t^\lambda(2 \y, 2\v)$ is even in $\lambda$ leads to the relation
$$
  W_t^+(\z,\w,i\eta) = W_t^-(\z,\w,-i\eta).
$$  

\par We refer to $W_t^+$ and $W_t^-$ as the {\it partial weight functions}. 
Their importance will become clear in the next subsection. 

\begin{rem} We will show in the appendix that both 
$W_t^+$ and $W_t^-$ attain positive and negative values. In addition 
we shall discuss their oscillatory behaviour. A more heuristic 
explanation of these phenomena might 
be the following: 
Both $W_t^+(i\y,i\v,i\eta)$ and $W_t^-(i\y,i\v,i\eta)$ 
satisfy the differential equation
\begin{equation} \label{eq:diff} 
 2 \frac{\partial}{\partial t} U = \left(\Delta + (1- |\y|^2 - |\v|^2)
   \frac{\partial^2}{\partial \eta^2} \right) U. 
\end{equation}
Indeed, this follows from a straightforward computation starting from
$$
\frac{\partial}{\partial t} p_t^\lambda (\y,\v) = 
\left(\Delta -  \frac{\lambda^2}{4}(|\y|^2 + |\v|^2) \right)
 p_t^\lambda(\y,\v) 
$$
for all $\lambda \ne 0$ (see \cite{T}). We note that the 
differential 
equation \eqref{eq:diff} is parabolic only for $|\y|^2 + |\v|^2 <1$. 
If $|\y|^2 + |\v|^2 >1$, then the right hand side of \eqref{eq:diff} 
resembles a wave equation which in turn might explain the oscillatory 
behaviour of $W_t^+$ and $W_t^-$ on the large scale.    
\end{rem}

\subsection{The image of the heat kernel transform}
\setcounter{equation}{0}

The objective of this section is to prove our main 
theorem: $\im \CH_t=\CB_t^+(\H_\C)\oplus \CB_t^-(\H_\C)$ is 
a sum of two weighted Bergman spaces. 

\par To exhibit the Bergman structure of the spaces $\CB_t^+(\H_\C)$ and 
$\CB_t^-(\H_\C)$  needs some preparation. 

\par First we define subspaces of $L^2(\H)$ by 
$$
  L^2_+(\H) = \{f \in L^2(\H): f^\lambda = 0, \quad \lambda \leq 0\}
$$
and 
$$
  L^2_-(\H) = \{f \in L^2(\H): f^\lambda = 0, \quad \lambda \geq 0\}.
$$
Notice that both subspaces are $\H$-invariant and 
$$
  L^2(\H)  = L^2_+(\H) \oplus L^2_-(\H).
$$ 

\par Next we recall some facts on the heat 
kernel transform on the real line. The heat kernel on $\R$ 
is given by 
$$q_t(x) =(4\pi t)^{-\frac{1}{2}} e^{-\frac{x^2}{4t}} \qquad (x\in \R)\, .$$
Define a weighted Bergman space on $\C$ by 
$$\CB_t(\C)=\{ g\in \CO(\C): \|g\|^2 =\int_\C |g(x+iy)|^2 e^{-\frac{y^2}{2t}}\, dx\, dy<\infty\}
$$
and recall that the mapping 
$$h_t: L^2(\R)\to \CB_t(\C), \ \ g\mapsto (f*q_t)^\sim$$
is (up to scale) an $\R$-equivariant isometric isomorphism. 

\par Set $\R^+=(0,\infty)$ and $\R^-=(-\infty,0)$.
With $L_+^2(\R)=\{ f\in L^2(\R):{\rm supp}\hat f\subseteq\R^+ \}$
and $L_-^2(\R)=\{ f\in L^2(\R):{\rm supp}\hat f\subseteq\R^- \}$  
we have $L^2(\R)=L_+^2(\R)\oplus L_-^2(\R)$. Finally, let 
us write $\CB_t^\pm(\C)=h_t(L_\pm^2(\R))$. Clearly we have 
$\CB_t(\C)=\CB_t^+(\C)\oplus \CB_t^-(\C)$. 

\par Let $R>0$. Denote by $B_R$ the open ball centered at $0$ with radius 
$R$ in $\C^n$. Further define $K_R=B_R\times B_R\times \C\subseteq\H_\C$ 
and note that $\bigcup_{R>0} K_R=\H_\C$. 
\par 
\par  We define $\CV_t^+(\H_\C)$ as the vector space consisting
of all holomorphic functions $F$ on $\H_\C$ such that 
\begin{itemize}
\item  $F|_{K_R}\in L^2 (K_R, |W_t^+|dz)$ for all $R>0$, 
\item  $\lim_{R\to \infty} \int_{K_R} |F(z)|^2 W_t^+(z) \ dz <\infty$, 
\item  $F(\z,\w, \cdot)\in \CB_t^+(\C)$ for all $\z,\w\in \C^n$. 
\end{itemize}

\par We endow $\CV_t^+(\H_\C)$ with a sesquilinear bracket 
\begin{equation}\label{bra}
\la F, G\ra_+ =\lim_{R\to \infty} \int_{K_R} F(z) \oline{G(z)} 
 W_t^+(z)\, dz \, , 
\end{equation}
for $F, G\in \CV_t^+(\H_\C)$.
Similarly one defines $\CV_t^-(\H_\C)$ and $\la\cdot,\cdot\ra_-$. 

\begin{rem} One might ask if one cannot define  $\CV_t^\pm (\H_\C)$ in a
simpler manner: avoid the exhaustion $\bigcup_{R>0}K_R=\H_\C$ and 
just require $|F|^2 W_t^\pm $  to be absolutely 
integrable on $\H_\C$. 
However, this will not work, and the reason for this is the bad oscillatory 
behaviour of $W_t^\pm$ (see the appendix). 
\end{rem}
A priori it is not clear that $\la F, F\ra_\pm \geq 0$. This 
will be shown next. 

\begin{lem} The bracket $\la\cdot, \cdot\ra_\pm $ induces 
on $\CV_t^\pm (\H_\C)$ a pre Hilbert space structure. 
\end{lem}

\begin{proof} It is sufficient to treat  the case  ``$+$'' only. 
All what is left to show is that $\la F,F\ra_+\geq 0$ and 
$\la F,F\ra_+=0$ if and only if $F=0$. 
\par Fix $F\in \CV_t^+(\H_\C)$. Then $F(\z,\w, \cdot)\in \CB_t^+(\C)$
implies the existence of a function $g(\z,\w,\cdot)\in L_+^2(\R)$ such that 
$$
F(\z,\w,\zeta)=h_t(g(\z,\w,\cdot))(\zeta)=\int_\R g(\z,\w,s) 
q_t(\zeta-s)\, ds\, .
$$
Therefore, up to an irrelevant constant only depending on $t$, 
the following equality holds: 
$$\int_\R F(\z,\w,\xi+i\eta) e^{i\lambda\xi}\, d\xi
=e^{\lambda\eta} e^{-t\lambda^2} g^\lambda(\z,\w)\ .$$
Consequently, as $W_t^+$ is independent of $\xi$, 
$$
\int_{K_R} |F(z)|^2 W_t^+(z)\, dz =\int_{B_R^2} \int_0^\infty\int_\R 
|g^\lambda(\z,\w)|^2  e^{2\lambda\eta} e^{-2t\lambda^2} W_t^+(\z,\w,i\eta)
\, d\eta\, d\lambda\, d\z\,d\w\ .
$$
In view of (\ref{1002}) we thus get 
$$\int_{K_R} |F(z)|^2 W_t^+(z)\, dz =\int_{B_R}\int_{B_R}\int_0^\infty 
|g^\lambda(\z,\w)|^2  W_t^\lambda(\z,\w) \, d\lambda\, d\z\,d\w\ .$$
But $W_t^\lambda\geq 0$ and so 
$$
\la F, F\ra_+=\lim_{R\to \infty} \int_{B_R}\int_{B_R} \int_0^\infty 
|g^\lambda(\z,\w)|^2  W_t^\lambda(\z,\w) \, d\lambda\, d\z\,d\w\geq 0
$$
and $\la F, F\ra_+=0$ if and only if $g^\lambda=0$ for all $\lambda$, i.e.
$F=0$. This completes the proof of the lemma.  
\end{proof}
Let us write $\CH_t^\pm$ 
for the heat kernel transform when restricted to $L_\pm^2(\H)$.
Define  Hilbert spaces of holomorphic functions by 
$\CB_t^\pm (\H_\C)=\im \CH_t^{\pm}$
and note that 
 
$$\im \CH_t=\CB_t^+(\H_\C)\oplus\CB_t^-(\H_\C)\, .$$ 
Let us remark that this decomposition can be also achieved using the 
Hilbert transform in the last variable.

\begin{thm} Let $t>0$. Then $\CB_t^\pm(\H_\C)$ is the Hilbert 
completion of $(\CV_t^{\pm}(\H_\C), \la \cdot, \cdot\ra_\pm)$
with $\la\cdot, \cdot\ra_\pm$ given by (\ref{bra}).   
\end{thm}

\begin{proof}  We restrict ourselves to the ``$+$''-case. 
Define a dense subspace of $L_+^2(\H)^0$ of $L_+^2(\H)$ by 

$$L_+^2(\H)^0=\{f\in L_+^2(\H): \lambda \mapsto f^\lambda\  
\ \hbox{compactly  supported in} \ (0,\infty)\}$$
We claim that $\CH_t(L_+^2(\H)^0)\subset \CV_t^+(\H)$. 
Let $f\in L_+^2(\H)^0$ and set $F=\CH_t^+(f)$.
Choose $a>0$ such that 
$f^\lambda=0 $ for $ \lambda$ outside of $(a^{-1}, a)$.
Proceeding as in Lemma 6.4 and using the estimate 
Proposition 6.3 (ii)  we see that  $ \CH_t^+(f)$ 
satisfies the first
condition in the definition of $\CV_t^+(\H_\C).$  Furthermore
(6.1.5) implies that
$$
\int_{K_R} |F(z)|^2 W_t^+(z)\, dz=\int_{B_R}\int_{B_R}\int_0^\infty
|H_t^\lambda(f^\lambda)(\z,\w)|^2 W_t^\lambda(\z,\w)\, d\lambda\, d\z\ d\w\,
 .$$
As $W_t^\lambda\geq 0$, it hence follows that 
$\int_{K_R} |F(z)|^2 W_t^+(z)\,
 dz$ is increasing in $R$. Similar reasoning as in (6.1.6) now shows that
$$\lim_{R\to \infty} \int_{K_R} |F(z)|^2 W_t^+(z)\, dz=\|f\|^2 <\infty\ .$$
Furthermore, for fixed $(\z,\w)$ we have $F(\z,\w,\cdot)\in \CB_t(\C)$
as a quick inspection of (6.1.5) shows. 
This proves our claim. 

\par As a byproduct of our reasoning above we have 
shown that $\CH_t^+: L_2^+(\H)^0\to \CV_t^+(\H)$ is an isometric
map. It remains to verify that each function 
$F\in \CV_t^+(\H_\C)$ can be written as 
$\CH_t^+(f)$ for some $f\in L_+^2(\H)$. 
Let  $g^\lambda(\z,\w)$ be the function
associated to $F$ as in the proof of Lemma 6.4. Then for almost all 
$\lambda$ there exists an $f^\lambda\in L^2(\R^{2n})$ such that 
$g^\lambda=H_t^\lambda(f^\lambda)$. It is easy to check that 
the prescription 
$$
f(\x,\u,\xi)=\int_{\R} e^{-i\lambda \xi} f^\lambda(\x,\u)\,  d\lambda
$$
defines a function in $L_+^2(\H)$ such that $\CH_t^+(f)=F$. 
This completes the proof of the theorem. 
\end{proof}

\section{Appendix: The oscillatory behaviour of the partial weight functions}
\setcounter{equation}{0}

This appendix is devoted to a closer study of the partial weight 
functions $W_t^\pm$. In particular we will detect ``good'' and ``bad''
directions for $W_t^\pm$, meaning rays in $H_\C$ on which $W_t^\pm$ stays 
positive resp. starts to oscillate. It is no loss 
of generality to treat the case of $W_t^+$ only. 

\par We start with an expicit formula for the function 
$W_t^+$.  Recall that the kernel $p_t^{\lambda}$ admits 
an expansion of the type \cite[p. 85]{T} 
$$
  p_t^{\lambda} (\y,\v) = (2\pi)^{-n} \lambda^n \sum_{k=0}^\infty
 e^{-(2k+n)|\lambda|t} L_k^{n-1}(\frac{|\lambda|}{2} (|\y|^2+|\v|^2)) 
   e^{-\frac{|\lambda|}{4}(|\y|^2+|\v|^2)}\, ,  
$$ 
where $L_k^{n-1}$ is the Laguerre polynomial of degree $k$ with parameter
$n-1$, which can be extended analytically to $\lambda + i s$ for 
$\lambda \ne 0$. Let $H_k(x)$ denote the Hermite polynomial, which can be 
defined by Rodrigue's formula $H_k(x) = (-1)^k e^{x^2} \frac{d^k}{d x^k} 
e^{-x^2}$. 

\begin{prop}
For $n =1$ and $\beta:=(y^2+v^2)$, 
\begin{align*}
& W_{t/2}^+(i y,i v,i\eta) = c \sqrt{\frac{\pi}{t}} \sum_{k=0}^\infty
   e^{- \frac{1}{4} \mu_k^2 } \\
& \quad \times \left[ \mu_k  \sum_{j=0}^k \frac{1}{j!}  
\left(\frac{\beta}{\sqrt{t}}\right)^j H_j (-\mu_k\sqrt{t}) \binom{k}{j} 
+ \frac{\beta}{t} \sum_{j=0}^{k-1} \frac{1}{j!}   
  \left(\frac{\beta}{\sqrt{t}}\right)^j
   H_j (-\mu_k\sqrt{t}) \binom{k}{j+1}\right] 
\end{align*}    
where $\mu_k =(2k+1 + (2\eta+\beta)/t)/2$. 
\end{prop}

\begin{proof}
The integral formula of $W_t^+$ shows that, for a fixed $\lambda>0$,  
\begin{align*}
& W_{t/2}^+(i y,iv, i\eta) 
 = c \int_\R  e^{t(\lambda+i s)^2} e^{-2\eta(\lambda+is)}  
   p_t^{\lambda+i s}(2y,2v) \,ds \\
& \qquad  = c\, e^{-(t + 2\eta)\lambda+t \lambda^2 - \lambda \beta} 
  \sum_{k=0}^\infty e^{-2k t \lambda} \int_\R  (\lambda+is) 
      e^{-t s^2} L_k^{n-1}(2 \lambda \beta + 2 i s \beta) \\
& \qquad   \hspace{2in} \times  e^{- i s t (2\lambda - 2k - \beta - 2\eta/t)}
       e^{- i s \beta} \,ds \\
& \qquad  =  c\, e^{-(t + 2\eta)\lambda+t \lambda^2 - \lambda \beta} 
  \sum_{k=0}^\infty e^{-2k t \lambda} (\lambda+\partial_\alpha) 
      L_k^{n-1}(2 \beta (\lambda +\partial_\alpha))
   \int_\R  e^{i \alpha s}e^{-t s^2} \,ds \\
& \qquad  =  c\,\sqrt{\frac{\pi}{t}}
  e^{-(t + 2\eta)\lambda+t \lambda^2 - \lambda \beta} 
  \sum_{k=0}^\infty e^{-2k t \lambda} (\lambda+\partial_\alpha) 
      L_k^{n-1}(2 \beta(\lambda + \partial_\alpha))
           e^{-\frac{1}{4t} \alpha^2} 
\end{align*}
where $\alpha = t( 2 \lambda - 2 k -1 - (2\eta + \beta)/t)$ and 
$\partial_\alpha = \partial /\partial_\alpha$. 

Using the Rodrigue's formula of the Hermite polynomials and 
the explicit formula of $L_k^\gamma$, we conclude that 
\begin{align*}
 L_k^{n-1}(2 \beta (\lambda+\partial_\alpha))
   e^{-\frac{1}{4t} \alpha^2} & = 
  \sum_{l=0}^k \frac{(-k)_l}{l! l!} (2\beta)^l (\lambda+ \partial_\alpha)^l
     e^{-\frac{1}{4t} \alpha^2} \\ 
&  = \sum_{l=0}^k \frac{(-k)_l}{l! l!} (2\beta)^l 
     \sum_{j = 0}^l \binom{l}{j} \frac{1}{(2\sqrt{t})^j} 
 (-1)^j e^{- \frac{1}{4t}\alpha^2} H_j \left(\frac{\alpha}{2\sqrt{t}}\right)\\
& = e^{- \frac{1}{4t}\alpha^2} \sum_{j=0}^k \frac{1}{j!} 
   \frac{\beta}{(\sqrt{t})^j} H_j \left(\frac{\alpha}{2\sqrt{t}}\right)
    L_{k-j}^j(2 \beta \lambda), 
\end{align*}  
upon changing summations, simplifying and using the explicit formula of 
$L_k^j$. Let $\alpha$ be fixed. It turns out that the generating function
of the above quantity is given by   
$$ 
 e^{- \frac{1}{4t}\alpha^2} \sum_{k=0}^\infty \left(
  \sum_{j=0}^k \frac{1}{j!} \frac{\beta}{(\sqrt{t})^j} 
      H_j \left(\frac{\alpha}{2\sqrt{t}}\right) 
     L_{k-j}^j(2 \beta \lambda)\right) s^k  
=\exp\left[{-\frac{2\beta s \mu}{1-s} - \frac{\beta s}{(1-s)\sqrt{t}} }\right]
$$     
where $\alpha  = 2t (\lambda - \mu)$. Since the generating function is 
independent of $\lambda$, this shows that the inner sum is in fact 
independent of $\lambda$. We can, in particular, set $\lambda = 0$ in the 
inner sum and set $\mu = (2k+1+ (2\eta + \beta)/t)/2$. Recall that 
$L_{k-j}^j(0) = \binom{k}{j}$. The change of variable from $\alpha$ to 
$\mu$ also leads to $\partial_\alpha = - \frac{1}{2t} \partial_\mu$. 
A simple computation then leads to 
\begin{align*}  
& (\lambda+\partial_\alpha) L_k^{n-1}(2 \beta (\lambda+\partial_\alpha))  
   e^{-\frac{1}{4t} \alpha^2}  
 = e^{-t(\lambda - \mu)^2} \\
& \quad \times \left[ \mu \sum_{j=0}^k \frac{1}{j!}  
   \left( \frac{\beta}{\sqrt{t}}\right)^j H_j (-\sqrt{t}\,\mu) \binom{k}{j}+ 
   \frac{\beta}{t} \sum_{j=0}^{k-1} \frac{1}{j!}  
   \left(\frac{\beta}{\sqrt{t}}\right)^j
   H_j (-\sqrt{t}\,\mu) \binom{k}{j+1}\right] 
\end{align*}      
from which the stated formula follows readily.
\end{proof}

We note that the formula proved above shows explicitly that $W_t^+$ is
independent of $\lambda$ without using the contour integral and Cauchy's
theorem. 

\begin{prop}
The function $W_t^+$ is positive in a neighborhood of $(0,0,0)$. Furthermore,
$W_t^+(0,0,i\eta)$ is non-negative for all $\eta$.
\end{prop}    

\begin{proof} 
Setting $\beta =0$ in the explicit formula of $W_t^+$ gives 
$$
W_{t/2}^+(0,0,i\eta) = c\sqrt{\frac{\pi}{t}} 
 \sum_{k=0}^\infty e^{-\frac{1}{4} (2k+1+2\eta/t)^2} 
  \left(2k+1+\frac{2\eta}{t}\right), 
$$
which is clearly positive if $\eta \ge 0$. Furthermore, if $\eta/t = -m$ for
$m\in \N$ then the sum can be written as 
$$
 \sum_{k=0}^\infty e^{-\frac{t}{4} (2k+1 - 2m)^2}(2k+1-2m) = 
 \sum_{k=0}^\infty e^{-\frac{t}{4} (2k+1)^2}(2k+1)  
  -  \sum_{k=1}^{m} e^{-\frac{t}{4} (2k-1)^2}(2k- 1)  
$$
which is strictly positive. Similarly, the sum is strictly positive if
$\eta/t=-m-1/2$. Hence, we are left with the case of $2 \eta /t = - 2m -1+r$,
where $0 < r < 1$. In this case, the sum becomes 
$$
 S_m:= 2\sum_{k=0}^\infty e^{-(k- m + r/2)^2 t }(k-m+r/2)  
     = 2 \sum_{k=0}^\infty e^{-(k + s)^2 t }(k+s)  -
       \sum_{k=1}^m e^{-(k- s)^2 t }(k-s)  
$$
where $0 < s=r/2 <1/2$. Set $g_k(s) = e^{-(k + s)^2 t }(k+s)  -
e^{-(k+1- s)^2 t }(k+1-s)$. It is easy to see that $g_k'(s) > 0$ for 
$0 < s < 1$. Hence $g_k$ is increasing. It follows that 
$$
  \sum_{k=0}^\infty e^{-\frac{t}{4} (2k+1)^2}(2k+1)  
  -  \sum_{k=1}^{\infty} e^{-\frac{t}{4} (2k-1)^2}(2k- 1)   
= \sum_{k=0}^\infty g_k(s) \ge \sum_{k=0}^\infty g_k(0) = 0,
$$
from which the stated result follows. 
 \end{proof}

However, the weight function $W_t^+(i y,i v,i \eta)$ is not non-negative 
for all $(y,v,\eta)$. In fact, if $2 \eta = - (y^2 + v^2)$, then 
$2 \eta + \beta =0$ and 
\begin{align*} 
& W_{t/2}^+(iy,iv,i\eta) = c \sqrt{\frac{\pi}{t}} \sum_{k=0}^\infty 
    e^{- (k+\frac{1}{2})^2} \left[ (k+\frac12) \sum_{j=0}^k \frac{1}{j!}   
\left(\frac{\beta}{\sqrt{t}}\right)^j H_j (-\sqrt{t}(k+\frac12)) \binom{k}{j}
 \right.  \\
& \hspace{2in} \left. + \frac{\beta}{t} \sum_{j=0}^{k-1} \frac{1}{j!}    
  \left(\frac{\beta}{\sqrt{t}}\right)^j 
   H_j (-\sqrt{t}(k+\frac12)) \binom{k}{j+1}\right].  
\end{align*}     
For each fixed $t$, this is a function of $\beta$ and it appears to be 
oscillatory. The graph for $t =1$ is shown below.

\smallskip
\input epsf.tex  
\centerline{\epsfxsize=4in \epsffile{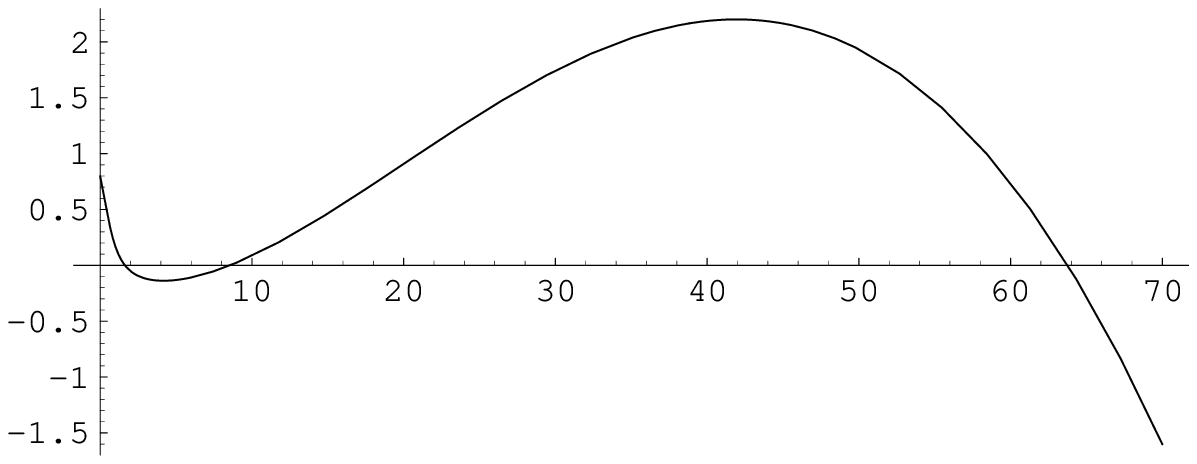} }
\smallskip

\noindent
The function oscillates in growing intervals and increasing amplitudes. 
To demonstrate the oscillatory nature of the function, what we have shown 
above is the function $W(iy,iv, -i\beta/2) /\log(2+\beta^2)$ without the
factor $c \sqrt{\pi}$. It is a function of $\beta$, where $\beta = y^2+v^2$.   


\end{document}